\documentclass[12pt]{amsart}
\usepackage[margin=2.0cm]{geometry}

\usepackage{amssymb}
\usepackage{amsfonts}
\usepackage{mathrsfs}
\usepackage{cite}
\usepackage{graphicx}

\usepackage{pgfplots}
\pgfplotsset{compat=1.15}
\usepackage{mathrsfs}
\usetikzlibrary{arrows}

\definecolor{ffqqqq}{rgb}{1,0,0}
\definecolor{qqzzqq}{rgb}{0,0.6,0}

\newcommand{\R}{{\mathbb R}}

\newtheorem{theorem}{Theorem}[section]

\newtheorem{lemma}[theorem]{Lemma}
\newtheorem{prop}[theorem]{Proposition}
\newtheorem{corollary}[theorem]{Corollary}

\newtheorem{example}[theorem]{Example}

\DeclareMathOperator{\arccosh}{\mathrm{arccosh}}
\DeclareMathOperator{\diam}{\mathrm{diam}}

\makeatletter
\DeclareFontFamily{U}{tipa}{}
\DeclareFontShape{U}{tipa}{m}{n}{<->tipa10}{}
\newcommand{\arc@char}{{\usefont{U}{tipa}{m}{n}\symbol{62}}}%

\newcommand{\arc}[1]{\mathpalette\arc@arc{#1}}

\newcommand{\arc@arc}[2]{%
  \sbox0{$\m@th#1#2$}%
  \vbox{
    \hbox{\resizebox{\wd0}{\height}{\arc@char}}
    \nointerlineskip
    \box0
  }%
}
\makeatother

\title[Extremal area of Reuleaux triangles]{Convex bodies of constant width in spaces of constant curvature and the extremal area of Reuleaux triangles}
\author{K\'aroly J. B\"or\"oczky, \'Ad\'am Sagmeister}
\address{Alfr\'ed R\'enyi Institute of Mathematics, Reltanoda u. 13-15, H-1053 Budapest, Hungary} \email{boroczky.karoly.j@renyi.hu}
\address{E\"otv\"os Lor\'and University, Institute of Mathematics, P\'azm\'any P\'eter s\'et\'any 1/c, Budapest, H-1117 Hungary} \email{sagmeister.adam@gmail.com }

\begin{document}

\maketitle

\begin{abstract}
Extending Blaschke and Lebesgue's classical result in the Euclidean plane, it has been recently proved in spherical and the hyperbolic cases, as well, that Reuleaux triangles have the minimal area among convex domains of constant width $D$. We prove an essentially optimal stability version of this statement in each of the three types of surfaces of constant curvature. In addition, we summarize the fundamental properties of convex bodies of constant width in spaces of constant curvature, and provide a characterization in the hyperbolic case in terms of horospheres.
\end{abstract}

\section{Introduction}

Let $\mathcal{M}^n$ be either the Euclidean space $\R^n$, hyperbolic space $H^n$ or spherical space $S^n$ for $n\geq 2$. We write $V_{\mathcal{M}^n}$ to denote the $n$-dimensional volume (Lebesgue measure) on $\mathcal{M}^n$,
and $d_{\mathcal{M}^n}(x,y)$ to denote the geodesic distance between $x,y\in\mathcal{M}^n$.
For $x,y\in \mathcal{M}^n$ where $x\neq - y$ if $\mathcal{M}^n=S^n$, we write $[x,y]_{\mathcal{M}^n}$ to denote the geodesic segment between $x$ and $y$ whose length is 
$d_{\mathcal{M}^n}(x,y)$. 
For $X\subset \mathcal{M}^n$ where $X$ is contained in an open hemisphere if $\mathcal{M}^n=S^n$, we say is $X$ convex if $[x,y]_{\mathcal{M}^n}\subset X$ for $x,y\in X$. In addition, $X$ is a convex body, if $X$ is compact and convex with non-empty interior.
For pairwise different points $x,y,z\in \mathcal{M}^n$ (which are pairwise not antipodal in the spherical case), let $\angle(x,y,z)$ be the angle of the geodesic segments $[x,y]_{\mathcal{M}^n}$ and $[y,z]_{\mathcal{M}^n}$ at $y$.
For $z\in\mathcal{M}^n$ and $r\geq 0$, we write
$B_{\mathcal{M}^n}(z,r)=\{x\in\mathcal{M}^n\colon d_{\mathcal{M}^n}(x,z)\leq r\}$ to denote the ball of radius $r$ and center $z$.

For a bounded set $X\subset \mathcal{M}^n$,  its diameter ${\rm diam}_{\mathcal{M}^n} X$ is the supremum of the geodesic distances $d_{\mathcal{M}^n}(x,y)$ for $x,y\in X$. We say that a bounded set $X\subset \mathcal{M}^n$ is complete if adding any point $z\not\in X$ to $X$, the diameter of $X\cup\{z\}$ is larger than ${\rm diam}_{\mathcal{M}^n}X$.
It is a well-known property of the Euclidean space (see Groemer \cite{Gro86}) that complete sets of diameter $D>0$ are the convex compact sets of constant width $D$; namely, the distance between any two parallel supporting hyperplanes is $D$.
Here we provide the analogue of the statement in the hyperbolic space. Hyperplanes in the Euclidean space can be thought as hypersurfaces orthogonal to a family of parallel lines. In the hyperbolic setting, hypersurfaces orthogonal to a family of parallel lines are horospheres (see Section~\ref{secconvex}).

\begin{theorem}\label{CompletenessHyptheo}
A set $K\subseteq H^n$ with 
${\rm diam}_{H^n}X=D>0$, is complete if and only if 
$K$ is compact convex and the distance between any two parallel supporting horospheres is $D$.
\end{theorem}

The two equivalent conditions in Theorem~\ref{CompletenessHyptheo} are also equivalent to Dekster's \cite{Dek95} definition of a convex body of constant width $D$ in a space of constant curvature (see Section~\ref{secconvex}).

 For $D>0$, let us consider the extremal values of the volume of a convex body of constant width $D$ in a space of constant curvature $\mathcal{M}^n$.
 According to the hyperbolic isodiametric inequality (see Schmidt \cite{Sch48,Sch49} and B\"or\"oczky, Sagmeister \cite{BoS20}) in
 $\mathcal{M}^n$, balls of diameter $D$ have the maximal volume among compact sets of diameter at most $D$. Therefore, balls have maximal volume among compact convex sets of constant width $D$.
 
 On the other hand, the minimal volume of a convex body of constant with $D$ (which exists by Blaschke's Selection Theorem~\ref{Blaschke} and Lemma~\ref{limit}) is not even known in the Euclidean space $\R^3$ (see Chakerian \cite{Cha66}, Schramm \cite{Sch88a}). Therefore, we focus on the minimal volume of a convex body  of constant width $D$ in $\mathcal{M}^2$.

A Reuleaux triangle in a surface of constant curvature is the intersection of three circular discs of radius $r>0$ whose centers are vertices of a regular triangle of side length $r$. It is a convex domain of constant width $r$ (see Section~\ref{secconvex} for definition and basic properties  of convex bodies of constant width).

The Blaschke--Lebesgue Theorem, due to
Blaschke \cite{Bla15} and 
Lebesgue \cite{Leb14}
states that amongst bodies of
constant width in the Euclidean plane, the Reuleaux triangle has the minimal area on the Euclidean plane (see Eggleston \cite{Egg52} for a particularly simple proof). The spherical version of the theorem was proved by Leichtweiss \cite{Lei05} based on some ideas of Blaschke. A new proof of the spherical case was recently published by K. Bezdek \cite{Bez21}. Here we prove the hyperbolic version of the Blaschke--Lebesgue Theorem. 

\begin{theorem}\label{BL-hyp}
Let $K\subseteq H^2$ be a body of constant width $D$ for some positive $D$, and let $U_D\subseteq H^2$ denote a Reuleaux triangle of width $D$. Then
$$
V_{H^2}\left(K\right)\geq V_{H^2}\left(U_D\right)
$$
with equality if and only if $K$ is congruent with $U_D$.
\end{theorem}

We note that Araujo \cite{Ara97} and
Leichtweiss \cite{Lei05} provided an involved arguments to verify Theorem~\ref{BL-hyp} if the boundary of $K$ is piecewise smooth. Our proof is much simpler and elementary.

For some related results in the hyperbolic space, see 
Alfonseca, Cordier, Florentin \cite{ACF21},
Gallego, Reventos, Solanes, Teufel \cite{GRST08},
Jerónimo-Castro, Jimenez-Lopez \cite{JCJL17}.

For related surveys on spherical convex bodies,
see Schramm \cite{Sch88b}, Lassak \cite{La20} and Lassak, Musielak\cite{LaM18}. In addition, Groemer \cite{Gro86} surveys properties of complete sets in Minkowski spaces.

The main result of our paper is the stability of the Blaschke--Lebesgue Inequality on any surface of constant curvature. Let $\mathcal{M}^2$ be either $\R^2$, $S^2$ or $H^2$. To measure the distance of two compact sets $X,Y\subset \mathcal{M}^2$, we use the Hausdorff distance
$$
\delta_H(X,Y)=\min\left\{r\geq 0\colon
B_{\mathcal{M}^2}(x,r)\cap Y\neq\emptyset\mbox{ and }
B_{\mathcal{M}^2}(y,r)\cap X\neq\emptyset
\mbox{ for any $x\in X$ and $y\in Y$}\right\}.
$$

\begin{theorem}\label{BL-stab}
Let $\mathcal{M}^2$ be either $\R^2$, $S^2$ or $H^2$ and let $D>0$ where $D<\frac{\pi}2$ if $\mathcal{M}^2=S^2$.
If $K\subset \mathcal{M}^2$ is a body of constant width $D$, $\varepsilon\geq 0$ and  
$$
V_{\mathcal{M}^2}\left(K\right)\leq (1+\varepsilon)V_{\mathcal{M}^2}\left(U_D\right),
$$
then there exists a Reuleaux triangle $U\subset \mathcal{M}^2$ of width $D$ such that
$\delta_H(K,U)\leq \theta\varepsilon$ where $\theta>0$ is an explicitly calculable constant depending on $D$ and $\mathcal{M}^2$.
\end{theorem}

We note that if $\mathcal{M}^2=\R^2$,
or $\mathcal{M}^2=H^2$ and $d\leq 1$, or
$\mathcal{M}^2=S^2$ and $d<\frac{\pi}2$, then
$\theta=c\cdot D$ can be chosen in Theorem~\ref{BL-stab} where $c>0$ is an absolute constant. However,
we have not made an attempt to actually calculate a constant $\theta$ in Theorem~\ref{BL-stab} because we can't expect one close formula for all $D>0$ in the hyperbolic case.
The reason is that if $\mathcal{M}^2=H^2$ and $R(U_D)$ is the circumradius of the Reuleaux triangle $U_D\subset H^2$ of width $D$, then
the Law of Sine yields $\sinh \frac{D}2/\sinh R(U_D)=\sin\frac{\pi}3$; therefore,
$$
\lim_{D\to 0}\frac{R(U_D)}{D/2}=\frac{2}{\sqrt{3}}
\mbox{ \ \ \ and \ \ \ }
\lim_{D\to \infty}\left(R(U_D)-\frac{D}{2}\right)
=\ln \frac{2}{\sqrt{3}}.
$$

Next we show that the order of the error estimate as a function of $\varepsilon$ is optimal in Theorem~\ref{BL-stab}.

\begin{example}\label{BL-stab-example}
Let $\mathcal{M}^2$ be either $\R^2$, $S^2$ or $H^2$ and let $D>0$ where $D<\frac{\pi}2$ if $\mathcal{M}^2=S^2$,
let $v_1,v_2,v_3\in \mathcal{M}^2$ be the vertices of a regular triangle with side length $D$, let
$p$ be the circumcenter of this triangle, and let $u_3$ be the point such that
$d_{\mathcal{M}^2}(v_3,u_3)=D$ and $p\in [v_3,u_3]_{\mathcal{M}^2}$ (see Section~\ref{secconstant} for notation). In particular, we may set
$U_D=\cap_{i=1,2,3} B_{\mathcal{M}^2}(v_i,D)$ to be our basic Reuleaux triangle.

For small $\eta>0$, let
$u_3^{(\eta)}$ be the point such that
$d_{\mathcal{M}^2}(u_3^{(\eta)},u_3)=\eta$ and $u_3\in [v_3,u_3^{(\eta)}]_{\mathcal{M}^2}$, and for $i=1,2$,
let $v_{3i}^{(\eta)}$ be the intersection of $\partial B_{\mathcal{M}^2}(u_3^{(\eta)},D)$ with the shorter arc
$\arc{v_3,v_i}$ of
$\partial B_{\mathcal{M}^2}(v_{3-i},D)$. We consider the
``Reuleaux pentagon''
$$
P^{(\eta)}=
B_{\mathcal{M}^2}(v_1,D)\cap
B_{\mathcal{M}^2}(v_2,D)\cap
B_{\mathcal{M}^2}(u_3^{(\eta)},D)\cap
B_{\mathcal{M}^2}(v_{31}^{(\eta)},D)\cap
B_{\mathcal{M}^2}(v_{32}^{(\eta)},D)
$$
that is a convex body of constant width $D$.

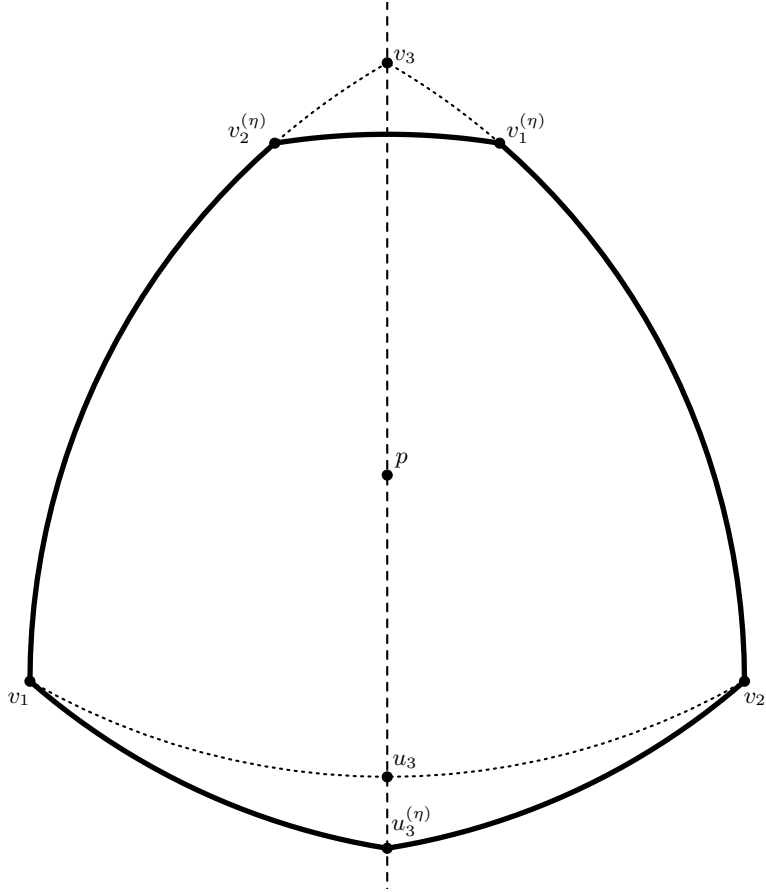
\begin{center}
\begin{tikzpicture}[line cap=round,line join=round,>=triangle 45,x=1cm,y=1cm, scale=9.5]
\clip(-0.072,-0.29) rectangle (1.144,0.95);
\draw [line width=.8pt,dashed] (0.5,-0.43494532702898564) -- (0.5,1.2223150329141554);
\draw [shift={(0.3426259992172425,0.7535644783924406)},line width=2pt]  plot[domain=4.870419932050983:5.429717650790183,variable=\t]({1*1*cos(\t r)+0*1*sin(\t r)},{0*1*cos(\t r)+1*1*sin(\t r)});
\draw [shift={(0.6573740007827575,0.7535644783924406)},line width=2pt]  plot[domain=3.995060309979196:4.554358028718395,variable=\t]({1*1*cos(\t r)+0*1*sin(\t r)},{0*1*cos(\t r)+1*1*sin(\t r)});
\draw [shift={(0,0)},line width=2pt]  plot[domain=0:0.853467656389403,variable=\t]({1*1*cos(\t r)+0*1*sin(\t r)},{0*1*cos(\t r)+1*1*sin(\t r)});
\draw [shift={(0.5,-0.23397459621556124)},line width=2pt]  plot[domain=1.412765375128603:1.7288272784611904,variable=\t]({1*1*cos(\t r)+0*1*sin(\t r)},{0*1*cos(\t r)+1*1*sin(\t r)});
\draw [shift={(1,0)},line width=2pt]  plot[domain=2.28812499720039:3.141592653589793,variable=\t]({1*1*cos(\t r)+0*1*sin(\t r)},{0*1*cos(\t r)+1*1*sin(\t r)});
\draw [shift={(0.5,0.8660254037844388)},line width=.8pt,dotted]  plot[domain=4.1887902047863905:5.235987755982988,variable=\t]({1*1*cos(\t r)+0*1*sin(\t r)},{0*1*cos(\t r)+1*1*sin(\t r)});
\draw [shift={(1,0)},line width=.8pt,dotted]  plot[domain=2.0943951023931953:2.28812499720039,variable=\t]({1*1*cos(\t r)+0*1*sin(\t r)},{0*1*cos(\t r)+1*1*sin(\t r)});
\draw [shift={(0,0)},line width=.8pt,dotted]  plot[domain=0.853467656389403:1.0471975511965976,variable=\t]({1*1*cos(\t r)+0*1*sin(\t r)},{0*1*cos(\t r)+1*1*sin(\t r)});
\begin{scriptsize}
\draw [fill=black] (0,0) circle (.2pt);
\draw[color=black] (-0.015,-0.025) node {$v_1$};
\draw [fill=black] (1,0) circle (.2pt);
\draw[color=black] (1.015,-0.025) node {$v_2$};
\draw [fill=black] (0.5,0.8660254037844388) circle (.2pt);
\draw[color=black] (0.525,0.875) node {$v_3$};
\draw [fill=black] (0.5,-0.23397459621556124) circle (.2pt);
\draw[color=black] (0.535,-0.196) node {$u_3^{\left(\eta\right)}$};
\draw [fill=black] (0.5,-0.13397459621556115) circle (.2pt);
\draw[color=black] (0.525,-0.112) node {$u_3$};
\draw [fill=black] (0.3426259992172425,0.7535644783924406) circle (.2pt);
\draw[color=black] (0.305,0.775) node {$v_2^{\left(\eta\right)}$};
\draw [fill=black] (0.6573740007827575,0.7535644783924406) circle (.2pt);
\draw[color=black] (0.695,0.775) node {$v_1^{\left(\eta\right)}$};
\draw [fill=black] (0.5,0.2886751345948129) circle (.2pt);
\draw[color=black] (0.52,0.31) node {$p$};
\end{scriptsize}
\end{tikzpicture}
\end{center}

Now there exist constants $\theta_1,\theta_2>0$ depending on $D$ and $\mathcal{M}^2$ such that
$V_{\mathcal{M}^2}\left(P^{(\eta)}\right)\leq (1+\theta_1\eta)V_{\mathcal{M}^2}\left(U_D\right)$,
and
$\delta_H(P^{(\eta)},U)\geq \theta_2\eta$ for any Reuleaux triangle $U\subset \mathcal{M}^2$ of width $D$.
Therefore, if 
$V_{\mathcal{M}^2}\left(P^{(\eta)}\right)= (1+\varepsilon)V_{\mathcal{M}^2}\left(U_D\right)$, 
then $\delta_H(P^{(\eta)},U)\geq \frac{\theta_2}{\theta_1}\varepsilon$ for any Reuleaux triangle $U\subset \mathcal{M}^2$ of width $D$ where
$\varepsilon$ can be arbitrarily small.
\end{example}

We note that stability versions of geometric inequalities is a well investigated subject. For example, essentially optimal stability version of the isoperimetric inequality in the terms of volume difference was proved by
Fusco, Maggi, Pratelli \cite{FMP08} in the Euclidean case (see also
Figalli, Maggi, Pratelli \cite{FMP09,FMP10}, and the best known constant is achieved in Kolesnikov, Milman \cite{KolMilsupernew}). The analogous result in 
the hyperbolic case is verified by B\"ogelein, Duzaar, Scheven \cite{BDS15},
and in the spherical case is proved by B\"ogelein, Duzaar, N. Fusco \cite{BDF17}. 
Stability versions of the isoperimetric inequality in $\R^n$ in the terms of Hausdorff distance have been verified by 
Diskant \cite{Dis73} and Groemer \cite{Gro88}
(see Groemer \cite{Gro93} for a survey on related stability results up to 1993).

\section{Spaces of constant curvature}
\label{secconstant}

Let $\mathcal{M}^n$ be either $\R^n$,  $H^n$ or $S^n$. Our focus is on the spherical- and hyperbolic space, and we assume that $S^n$ is embedded into $\R^{n+1}$ the standard way,
and $H^n$ is embedded into $\R^{n+1}$ using the hyperboloid model. We write $\langle \cdot,\cdot\rangle$ to denote
the standard scalar product in $\R^{n+1}$, 
and write 
$z^\bot=\{x\in\R^{n+1}\colon\langle x,z\rangle=0\}$ for a $z\in\R^{n+1}\setminus\left\{o\right\}$.
Fix an $e\in S^n$. In particular, we have
\begin{eqnarray*}
S^n&=& \{x\in\R^{n+1}\colon\langle x,x\rangle=1\}\\
H^n&=& \{x+te\colon x\in e^\bot \mbox{ and }t\geq 1\mbox{ and }t^2-\langle x,x\rangle=1\}.
\end{eqnarray*}

For $H^n$, we also consider the following symmetric bilinear form $\mathcal{B}$ on $\R^{n+1}$: If $x=x_0+te\in\R^{n+1}$ and $y=y_0+se\in\R^{n+1}$ for
$x_0,y_0\in e^\bot$ and $t,s\in\R$, then
$$
\mathcal{B}(x,y)=ts-\langle x_0,y_0\rangle.
$$
In particular,
\begin{equation}
\label{Hn}
\mathcal{B}(x,x)=1\mbox{ for $x\in H^n$}.
\end{equation}
We note that the isometries of $\R^n$ are the maps of the form $x\mapsto Ax+b$ where $A\in O(n)$ and $b\in \R$,
the isometries of $S^n\subset \R^{n+1}$ are the maps of the form $x\mapsto Ax$ where $A\in O(n+1)$,
and the  isometries of $H^n\subset \R^{n+1}$ are the maps of the form $x\mapsto Ax$ where $A\in {\rm GL}(n+1,\R)$
leaves $\mathcal{B}(\cdot,\cdot)$ invariant and $\langle Ae,e\rangle>0$. The isometry group of each
$\mathcal{M}^n$ acts transitively, and the subgroup fixing a $z\in \mathcal{M}^n$ is isomorphic to $O(n-1)$.

Again let $\mathcal{M}^n$ be either $\R^n$,  $H^n$ or $S^n$ using the models as above for $H^n$ and $S^n$. For
$z\in \mathcal{M}^n$, we define the tangent space $T_z$ as
\begin{eqnarray*}
T_z&=&\{x\in\R^{n+1}\colon\mathcal{B}(x,z)=0\} \mbox{ if $\mathcal{M}^n=H^n$}\\
T_z&=&z^\bot\subset \R^{n+1} \mbox{ if $\mathcal{M}^n=S^n$}\\
T_z&=&\R^n \mbox{ if $\mathcal{M}^n=\R^n$}.
\end{eqnarray*}  
We observe that $T_z$ is an $n$-dimensional real vector space equipped with the scalar product
$-\mathcal{B}(\cdot,\cdot)$ if $\mathcal{M}^n=H^n$, and with the scalar product 
$\langle \cdot,\cdot\rangle$ if $\mathcal{M}^n=S^n$ or $\mathcal{M}^n=\R^n$.

For $z\in \mathcal{M}^n$ and unit vector $u\in T_z$, the geodesic line $\ell$ passing through $z$ and determined by $u$ consists of the points
\begin{equation}
\label{geodesic}
p_t=\left\{
\begin{array}{rl}
z\,{\rm cosh}\,t+u\,{\rm sinh}\,t& \mbox{ \ if $\mathcal{M}^n=H^n$}\\
z\cos t+u\sin t& \mbox{ \ if $\mathcal{M}^n=S^n$}\\
z+tu& \mbox{ \ if $\mathcal{M}^n=\R^n$}
\end{array}\right.
\end{equation}
for $t\in \R$. Here the map $t\mapsto p_t$ is bijective onto $\ell$ and satisfies 
$d_{\mathcal{M}^n}(z,p_t)=|t|$ for $t\in \R$ if $\mathcal{M}^n=H^n$ or $\mathcal{M}^n=\R^n$,
and for $t\in(-\pi,\pi]$ if $\mathcal{M}^n=S^n$. If $t>0$  provided $\mathcal{M}^n=H^n$ or $\mathcal{M}^n=\R^n$,
or $0<t<\pi$ provided $\mathcal{M}^n=S^n$, then we say that $u$ points towards $p_t$ along the geodesic segment
$$
[z,p_t]_{\mathcal{M}^n}=\{p_s\colon 0\leq s\leq t\}
$$
of length $t$.

A hyperplane $H$ in $\mathcal{M}^n$ passing through the point $z\in \mathcal{M}^n$ and having exterior unit normal
$u\in T_z$, and
 the corresponding half-spaces $H^-$ and $H^+$ are defined as follows: $H^-=\mathcal{M}^n\setminus{\rm int} H^+$,
\begin{equation}
\label{hyperplane}
\begin{array}{lll}
H=\{x\in H^n\colon\mathcal{B}(x,u)=0\}& \mbox{ }H^+=\{x\in H^n\colon-\mathcal{B}(x,u)\geq 0\} &\mbox{ \ if $\mathcal{M}^n=H^n$}\\
H=\{x\in S^n\colon\langle x,u\rangle=0\}& \mbox{ }H^+=\{x\in S^n\colon\langle x,u\rangle\geq 0\}&\mbox{ \ if $\mathcal{M}^n=S^n$}\\
H=\{x\in \R^n\colon\langle x,u\rangle=\langle z,u\rangle\}& 
\mbox{ }H^+=\{x\in \R^n\colon\langle x,u\rangle\geq \langle z,u\rangle\} &\mbox{ \ if $\mathcal{M}^n=\R^n$}.
\end{array}
\end{equation}
The reflection $\sigma_H$ through $H$ is the unique isometry of $\mathcal{M}^n$ different from the identity fixing the points of $H$. In particular, for $x\in \mathcal{M}^n$ where $x\neq \pm u$ in the case $\mathcal{M}^n=S^n$, $H$ is the hyperplane perpendicularly bisecting the segment $[x,\sigma_Hx]_{\mathcal{M}^n}$ 
($H$ going through the midpoint of $[x,\sigma_Hx]_{\mathcal{M}^n}$).

An important tool to obtain convex bodies with extremal properties is the Blaschke Selection Theorem. 
First we impose a metric on compact subsets. Let $\mathcal{M}^n$ be either $\R^n$,  $H^n$ or $S^n$.
For a compact set $C\subset\mathcal{M}^n$ and $z\in \mathcal{M}^n$, 
we set $d_{\mathcal{M}^n}(z,C)=\min_{x\in C}d_{\mathcal{M}^n}(z,x)$. Depending on the space this geodesic distance is either
$$
d_{H^n}\left(x;y\right)=\arccosh\left(\mathcal{B}\left(x;y\right)\right)
$$
or
$$
d_{S^n}\left(x;y\right)=2\cdot\arcsin\left(\frac{d_{\mathbb{R}^{n+1}}\left(x;y\right)}{2}\right),
$$
and in $\mathbb{R}^n$ we use the Euclidean metric. For any non-empty compact set $C_1,C_2\subset\mathcal{M}^n$, we define their Hausdorff distance
$$
\delta_{\mathcal{M}^n}(C_1,C_2)=\max\left\{
\max_{x\in C_2}d_{\mathcal{M}^n}(x,C_1),
\max_{y\in C_1}d_{\mathcal{M}^n}(y,C_2)\right\}.
$$
The Hausdorff distance is a metric on the space of compact subsets in $\mathcal{M}^n$.
We say that a sequence $\{C_m\}$ of compact subsets of $\mathcal{M}^n$ is bounded if there is a ball containing every $C_m$.
For compact sets $C_m,C\subset \mathcal{M}^n$, we write $C_m\to C$ to denote if
the sequence $\{C_m\}$ tends to $C$ in terms of the Hausdorff distance.

For any $z\in \mathcal{M}^n$ and $r>0$, let
$$
B_{\mathcal{M}^n}(z,r)=\{x\in \mathcal{M}^n\colon d_{\mathcal{M}^n}(x,z)\leq r\}
$$ 
be the $n$-dimensional ball centered at $z$ where it is natural to assume $r<\pi$ if $\mathcal{M}^n=S^n$. When it is clear from the context what space we consider, we drop the subscript referring to the ambient space. In order to speak about the volume of a ball of radius $r$, we fix a reference point $z_0\in \mathcal{M}^n$ where $z_0=o$ the origin if $\mathcal{M}^n=\R^n$. The volume of the unit ball $B_{\R^n}\left(o,1\right)$ is denoted by $\kappa_n$.

The following is well-known  (see B\"or\"oczky, Sagmeister \cite{BoS20}).

\begin{lemma}
\label{Hausdorffconvergence}
For compact sets $C_m,C\subset \mathcal{M}^n$ where $\mathcal{M}^n$ is either $\R^n$,  $H^n$ or $S^n$, we have $C_m\to C$ if and only if
\begin{description}
\item[(i)] assuming $x_m\in C_m$, the sequence $\{x_m\}$ is bounded and any accumulation point of  $\{x_m\}$  
lies in $C$;
\item[(ii)] for any $y\in C$, there exist $x_m\in C_m$ for each $m$ such that $\lim_{m\to \infty}x_m=y$.
\end{description}
\end{lemma}
The space of compact subsets of $\mathcal{M}^n$ is locally compact according to the Blaschke Selection Theorem
(see R. Schneider \cite{Sch14}).

\begin{theorem}[Blaschke]
\label{Blaschke}
If $\mathcal{M}^n$ is either $\R^n$,  $H^n$ or $S^n$,
then any bounded sequence of compact subsets of $\mathcal{M}^n$ has
a convergent subsequence.
\end{theorem}

For convergent sequences of compact subsets of $\mathcal{M}^n$, we have the following 
(see B\"or\"oczky, Sagmeister \cite{BoS20}).

\begin{lemma}
\label{limit}
Let $\mathcal{M}^n$ be either $\R^n$,  $H^n$ or $S^n$, and
let the sequence $\{C_m\}$ of compact subsets of $\mathcal{M}^n$ tend to $C$.
\begin{description}
\item[(i)] ${\rm diam}_{\mathcal{M}^n}\,C=\lim_{m\to \infty}{\rm diam}_{\mathcal{M}^n}\,C_m$
\item[(ii)] $V_{\mathcal{M}^n}(C)\geq \limsup_{m\to \infty} V_{\mathcal{M}^n}(C_m)$
\end{description}
\end{lemma}

Recall that for any compact set  $X\subset \mathcal{M}^n$ and $\varrho\geq 0$, the
parallel domain is
$$
X^{(\varrho)}=\{z\in \mathcal{M}^n\colon\exists x\in X\mbox{ \ with \ }
d_{\mathcal{M}^n}(x,z)\leq \varrho\}=\bigcup \{B(x,\varrho)\colon x\in X\}.
$$
The triangle inequality and considering $x,y\in X$ with 
$d_{\mathcal{M}^n}(x.y)={\rm diam}_{\mathcal{M}^n}\, X$ show that
\begin{equation}
\label{parallel-diam}
{\rm diam}_{\mathcal{M}^n}\, X^{(\varrho)}=2\varrho+{\rm diam}_{\mathcal{M}^n}\, X.
\end{equation}
We discuss parallel sets based on Benyamini \cite{Ben84}.

\begin{lemma}
\label{limit-parallel}
If $\varrho\geq 0$ and  $\mathcal{M}^n$ is either $\R^n$,  $H^n$ or $S^n$, and
the sequence $\{C_m\}$ of compact subsets of $\mathcal{M}^n$ tends to $C$, then
\begin{description}
\item[(i)] $\left\{C_m^{(\varrho)}\right\}$ tends to $C^{(\varrho)}$;
\item[(ii)] ${\rm diam}_{\mathcal{M}^n}\,C^{(\varrho)}=
\lim_{m\to \infty}{\rm diam}_{\mathcal{M}^n}\,C_m^{(\varrho)}$;
\item[(iii)] $V_{\mathcal{M}^n}\left(C^{(\varrho)}\right)\geq \limsup_{m\to \infty} V_{\mathcal{M}^n}\big(C_m^{(\varrho)}\big)$;
\item[(iv)] for any $\varepsilon>0$,
$V_{\mathcal{M}^n}\left(C^{(\varrho)}\right)\leq \liminf_{m\to \infty} V_{\mathcal{M}^n}\big(C_m^{(\varrho+\varepsilon)}\big)$.
\end{description}
\end{lemma}
\proof We deduce (i) from Lemma~\ref{Hausdorffconvergence}, and in turn (ii)
from (\ref{parallel-diam}) and (iii)
 from 
Lemma~\ref{limit} (ii).

For (iv), we only observe that $C^{(\varrho)}\subset C_m^{(\varrho+\varepsilon)}$ if $m$ is large.
\endproof


\section{Inradius and circumradius of a convex body}

For a compact convex set $K$ in $\mathcal{M}^n$ 
where $\mathcal{M}^n$ is either $\R^n$,  $H^n$ or $S^n$, we  define the inradius and circumradius of $K$ as the maximal (and minimal) radius of a ball contained in (or containing) $K$. We use the notations $r\left(K\right)$ and $R\left(K\right)$, so
$$
r\left(K\right)=\max\left\{r\geq 0\colon\exists z\in\mathcal{M}^n\colon B_{\mathcal{M}^n}\left(z,r\right)\subseteq K\right\}
$$
and
$$
R\left(K\right)=\min\left\{R\geq 0\colon\exists z\in\mathcal{M}^n\colon K\subseteq B_{\mathcal{M}^n}\left(z,r\right)\right\}.
$$

We fix a regular $d$-dimensional simplex in $\mathcal{M}^n$ and of edge length $D$ is denoted by $\Delta_{\mathcal{M}^n}\left(D\right)$. Note that $\diam_{\mathcal{M}^n}\left(\Delta_{\mathcal{M}^n}\left(D\right)\right)=D$.

The following  version of Jung's Theorem is well known (see say Dekster \cite{Dek95J} in the hyperbolic case), but we verify it to have the exact statements we need later.

\begin{lemma}\label{jung}
If $K\subset \mathcal{M}^n$ a compact set of diameter $D>0$
where $\mathcal{M}^n$ is either $\R^n$,  $H^n$ or $S^n$
and $D<\frac{\pi}2$ in the case $\mathcal{M}^n=S^n$, then
\begin{description}
    \item[(i)] there is a unique point $p$ such that $K\subset B_{\mathcal{M}^n}\left(p,R\left(K\right)\right)$;
    \item[(ii)] for $2\leq k\leq n+1$, there is a set $\left\{q_1;\ldots;q_k\right\}\subset\partial_{\mathcal{M}^n}K\cap\partial_{\mathcal{M}^n}B_{\mathcal{M}^n}\left(p,R\left(K\right)\right)$ such that $p\in\mathrm{relint}\left[q_1,\ldots,q_k\right]_{\mathcal{M}^n}$ where $\left[q_1,\ldots,q_k\right]_{\mathcal{M}^n}$ is a $\left(k-1\right)$-dimensional simplex in $\mathcal{M}^n$;
    \item[(iii)] $$\frac{D}{2}\leq R\left(K\right)\leq R\left(\Delta_{\mathcal{M}^n}\left(D\right)\right)<D;$$
    \item[(iv)]  $R\left(K\right)=R\left(\Delta_{\mathcal{M}^n}\left(D\right)\right)$ holds if and only if $k=n+1$  and $\left[q_1,\ldots,q_{n+1}\right]_{\mathcal{M}^n}$ is congruent with $\Delta_{\mathcal{M}^n}\left(D\right)$ in (ii).
\end{description}
\end{lemma}
\proof To simplify notation, we drop the reference to $\mathcal{M}^n$ in the formulas.

Since the intersection of two different balls of radius $\varrho>0$ is contained in a ball of radius less than $\varrho$, we deduce that the circumcenter $p$ of $K$ is unique. 

For $\Xi=\partial K\cap\partial B(p,R(K))$, we claim that
\begin{equation}
\label{XiRK}
p\in {\rm conv}\,\Xi.
\end{equation}
We suppose that \eqref{XiRK} does not hold, and hence there exists a closed half-space $H^-$ containing $\Xi$ with $p\not\in H^-$, and let $H^+$  be the complementary closed half-space. It follows that $H^+\cap K\subset{\rm int}\,B(p,R(K))$; therefore, there exists a $p'$ closer to $H^-$ than $p$ such that $K\subset{\rm int}\,B(p',R(K))$. This contradiction with the definition of $R(K)$ yields \eqref{XiRK}. 

It follows from \eqref{XiRK} and Carath\'eodory's theorem 
that there exist
$2\leq k\leq n+1$ and $q_1;\ldots;q_k\in\partial K\cap\partial B_{\mathcal{M}^n}\left(p,R\left(K\right)\right)$ such that $p\in\left[q_1,\ldots,q_k\right]$. Assuming that $k$ is the possible smallest number with this property, we deduce that
$\left[q_1,\ldots,q_k\right]$ is a $\left(k-1\right)$-dimensional simplex and
$p\in\mathrm{relint}\left[q_1,\ldots,q_k\right]$.

Now readily $R(K)\geq D/2$. If $k=2$, then $R(K)=D/2$. Therefore, to verify (iii) and (iv), we may assume that $k\geq 3$.

Let
$u_1,\ldots,u_k\in T_p$ the unit vectors pointing towards
$q_1,\ldots,q_k$, respectively. 
It follows from 
\eqref{XiRK} that the origin $o$ of $T_p$ satisfies that
$o\in \left[u_1,\ldots,u_k\right]_{T_p}$. Since 
$u_1,\ldots,u_k$ are contained in a $(k-2)$-sphere in $T_p$ by $p\in\mathrm{relint}\left[q_1,\ldots,q_k\right]$, 
\begin{equation}
\label{uiujo}
\max_{i,j=1,\ldots,k} \angle(u_i,o,u_j) \geq\arccos \frac{-1}{k-1}   
\end{equation}
according to the fundamental result about packing $k$ equal spherical balls on a $\left(k-2\right)$-sphere (see Theorem~6.1.1 in B\"or\"oczky \cite{Bor04}). In addition, equality holds in \eqref{uiujo} if and only if
$u_1,\ldots,u_k$ are vertices of a regular $(k-1)$-simplex in $T_p$. Since $\angle(u_i,o,u_j)=\angle(q_i,p,q_j)$, we conclude (iii) and (iv) by the Law of Cosines for sides in $\mathcal{M}^2$.
\hfill $\Box$\\

Concerning the inradius, we have the following result.

\begin{lemma}\label{closest-points-on-the-boundary}
If $K\subseteq\mathcal{M}^n$ is a convex body 
where $\mathcal{M}^n$ is either $\R^n$,  $H^n$ or $S^n$, and $B_{\mathcal{M}^n}\left(w,r\left(K\right)\right)\subset K$, then for $2\leq k\leq n+1$ there are points $t_1,\ldots,t_k\in\partial_{\mathcal{M}^n}K\cap\partial_{\mathcal{M}^n}B_{\mathcal{M}^n}\left(w,r\left(K\right)\right)$ such that $\left[t_1,\ldots,t_k\right]_{\mathcal{M}^n}$ is a $\left(k-1\right)$-dimensional simplex and $w\in{\rm{relint}}\left[t_1,\ldots,t_k\right]_{\mathcal{M}^n}$.
\end{lemma}

\proof To simplify notation, we drop the reference to $\mathcal{M}^n$ in the formulas.

For $\Xi=\partial K\cap\partial B(w,r(K))$, we claim that
\begin{equation}
\label{XiinradiusK}
w\in {\rm conv}\,\Xi.
\end{equation}
We suppose that \eqref{XiinradiusK} does not hold, and hence there exists a closed half-space $H^-$ containing $\Xi$ with $w\not\in H^-$, and let $H^+$  be the complementary closed half-space. It follows that $H^+\cap B(p,r(K))\subset{\rm int}\,K$; therefore, there exists a $w'$ farther from $H^-$ than $w$ such that $B(w',R(K))\subset{\rm int}\,K$. This contradiction with the definition of $r(K)$ yields \eqref{XiinradiusK}. 

It follows from \eqref{XiinradiusK} and Carath\'eodory's theorem 
that there exist
$2\leq k\leq n+1$ and $t_1;\ldots;t_k\in\partial K\cap\partial B_{\mathcal{M}^n}\left(w,r\left(K\right)\right)$ such that $w\in\left[t_1,\ldots,t_k\right]$. Assuming that $k$ is the possible smallest number with this property, we deduce that
$\left[t_1,\ldots,t_k\right]$ is a $\left(k-1\right)$-dimensional simplex and
$w\in\mathrm{relint}\left[t_1,\ldots,t_k\right]$.
\hfill $\Box$

\noindent {\bf Remark} Note that contrary to the circumscribed ball, the inscribed ball of a compact convex set $K$ is not necessarily unique in the Euclidean or the hyperbolic space. For example, take a segment $s$ and let $K$ be the set of points of distance at most $r$ from $s$ for a fixed $r>0$.

\section{Complete sets in a space of constant curvature}
\label{secconvex}

Let $\mathcal{M}^n$ be either $\R^n$,  $H^n$ or $S^n$.
We call $X\subset \mathcal{M}^n$ convex if $[x,y]_{\mathcal{M}^n}\subset X$
for any $x,y\in X$, and in addition, we also assume that $X$ is contained in an open hemisphere if $\mathcal{M}^n=S^n$.
 For $Z\subset \mathcal{M}^n$ where we assume that 
$Z$ is contained in an open hemisphere if $\mathcal{M}^n=S^n$, the convex hull 
${\rm conv}_{\mathcal{M}^n} Z$ is the intersection of all convex sets containing $Z$. Similarly to the notation of geodesic segments which are just the convex hulls of their endpoints, for a finite set of points $\left\{q_1,\ldots,q_k\right\}$ in $\mathcal{M}^n$ we use the notation $\left[q_1,\ldots,q_k\right]_{\mathcal{M}^n}$ for the convex hull of the set.

 We remark that $B_{S^n}(z,r)$ is  convex if $r\in\left(0,\frac{\pi}2\right)$, and
 $B_{S^n}(z,r)$ is not convex if $r\in\left[\frac{\pi}2,\pi\right)$.

We call a $K\subseteq\mathcal{M}^n$ set \emph{complete} if for each $y\in\mathcal{M}^n\setminus K$, $\diam_{\mathcal{M}^n}\left(K\cup\left\{y\right\}\right)>\diam_{\mathcal{M}^n}\left(K\right)$. Readily, any complete set is closed. A complete set $\widetilde{K}$ is called a \emph{completion} of $K$ if $K\subseteq\widetilde{K}$ and they have the same diameter.
The Zorn lemma yields the existence of completions.

\begin{lemma}\label{Zorn}
Assuming the Axiom of Choice,
if $X\subset \mathcal{M}^n$ is a set of diameter at most $D>0$
where $\mathcal{M}^n$ is either $\R^n$,  $H^n$ or $S^n$
and $D<\frac{\pi}2$ in the case $\mathcal{M}^n=S^n$, then
there exists a complete set $\widetilde{X}\supset X$
of diameter $D$.
\end{lemma}

Another useful statement about completions is the following:

\begin{lemma}\label{diametercapballs}
If $X\subset \mathcal{M}^n$ is a set of diameter at most $D>0$
where $\mathcal{M}^n$ is either $\R^n$,  $H^n$ or $S^n$
and $D<\frac{\pi}2$ in the case $\mathcal{M}^n=S^n$, then
$$
{\rm diam}_{\mathcal{M}^n} Y\leq D
\mbox{ \ and $X\subset Y$ for }Y=\bigcap\left\{B_{\mathcal{M}^n}(z,D)\colon X\subset B_{\mathcal{M}^n}(z,D)\right\}.
$$
\end{lemma}
\proof Let $y,z\in Y$. Since $y\in B_{\mathcal{M}^n}(x,D)$ for any $x\in X$, we deduce that $X\subset B_{\mathcal{M}^n}(y,D)$, and hence $z\in Y\subset B_{\mathcal{M}^n}(y,D)$.
\hfill $\Box$\\

For any closed convex set $K\subset \mathcal{M}^n$ where $K$ is contained in an open hemisphere if $\mathcal{M}^n=S^n$,
and for any $z\in\partial K$, there exists a supporting half-space $H^+\subset \mathcal{M}^n$ such that $K\subset H^+$ and
$z\in\partial H^+$. In this case, any non-zero vector $v\in T_z$ orthogonal to $\partial H$ and exterior normal to $H^+$ is called an exterior normal to $K$ at $z$. The convexity and closedness of $K$ yields that the set $N_{\mathcal{M}^n}(K,z)\subset T_z$ consisting of the origin of $T_z$ and any exterior normal to $K$ at $z$ is a closed convex cone.

According to Dekster \cite{Dek95}, a convex body $K$ in $\mathcal{M}^n$ is of \emph{constant width} $D$ where $D>0$ (and in the spherical case we also assume $D<\frac{\pi}{2}$) if $\diam_{\mathcal{M}^n}\left(K\right)=D$, and for any $p\in\partial_{\mathcal{M}^n}K$ and any outer unit normal $v\in T_p$ of $K$ in $p$ there is a geodesic segment $\left[p,q\right]_{\mathcal{M}^n}\subseteq K$ along $-v$ of length $D$.

\begin{prop}
\label{completeness}
Let $\mathcal{M}^n$ be either $\R^n$, $S^n$ or $H^n$, and let
$K\subset\mathcal{M}^n$ be a compact set of diameter $D>0$ where $D<\frac{\pi}{2}$ if $\mathcal{M}^n=S^n$. The following are equivalent:
\begin{description}
\item[(i)] $K$ is of constant width $D$;
        \item[(ii)] $K$ is complete;
    \item[(iii)] $K=\bigcap_{x\in K}{B_{\mathcal{M}^n}\left(x,D\right)}$.
\end{description}
\end{prop}
\proof To simplify notation, we drop the reference to $\mathcal{M}^n$ in the formulas.

To show that (i) yields (ii), we assume that
$K$ is of constant width $D$ and $z\not\in K$. There exists a closest point $p\in K$ to $z$. Let $v\in T_p$ be the unit vector pointing along $[p,z]$, and hence $v$ is an exterior unit normal to $K$ at $p$. It follows that there is a $q\in K$ with $d(p,q)=D$ such that $-v$ points towards $q$ and $[p,q]\subset K$. Therefore, ${\rm diam} (D\cup\{z\})>D$, proving (ii). 

Next if $K$ is complete, then (iii) directly follows from
Lemma~\ref{diametercapballs}.

Finally, we prove that (iii) implies (i). First of all, (iii) yields that $K$ is a convex body. Let $p\in \partial K$. It follows that there exists some $x_0\in K$ such that
$K\subset B(x_0,D)$ and $p\in\partial B(x_0,D)$. Let
$\widetilde{N}(K,p)\subset T_p$ consist of the origin $o$ of $T_p$ and any exterior normal at $p$ to any ball $B(x,D)$ where $K\subset B(x,D)$, $x\in K$ and $p\in\partial B(x,D)$. 
Since $K$ is compact and ${\rm int}\,K\neq \emptyset$, it follows that
$\widetilde{N}(K,p)$ is a closed convex cone such that $o$ is an apex ($\widetilde{N}(K,p)\setminus\{o\}$ is contained in an open half-space of $T_p$ such that $o$ lies on the boundary). We claim that
\begin{equation}
\label{NtildeN}
\widetilde{N}(K,p)=N(K,p).
\end{equation}
We suppose that \eqref{NtildeN} does not hold, and seek a contradiction.
Since readily $\widetilde{N}(K,p)\subset N(K,p)$, there exists a non-zero exterior normal $v$ to $K$ at $p$ such that
$v\not\in \widetilde{N}(K,p)$. As $o$ is an apex of $\widetilde{N}(K,p)$, $v$ can be strictly separated from $\widetilde{N}(K,p)$, namely, there exists some non-zero $u\in T_p$ such that 
\begin{equation}
\label{uNtildeN}
\angle(u,o,v)<\frac{\pi}2\mbox{ \ and \ }
\angle(u,o,w)>\frac{\pi}2 \mbox{ \ for any \ } w\in \widetilde{N}(K,p)\setminus\{o\}.
\end{equation}

We choose a $q\not\in K$ such that $u$ points towards $q$ along $[p,q]$ (and $d(p,q)<\frac{\pi}2-D$ if $\mathcal{M}^n=S^n$).
As $\angle(u,o,v)<\frac{\pi}2$, we have $[p,q]\cap K=\{p\}$.
We consider a sequence $y_n\in [p,q]\setminus\left\{p\right\}$ tending to $p$, and hence $y_n\not\in K$. It follows from (iii) that for each $y_n$ there exists $x_n\in K$ with
$y_n\not \in B(x_n,D)$. We may assume that $x_n$ tends to $\widetilde{x}_0\in K$ where 
$d(p,\widetilde{x}_0)=D$ and
$[p,q]\cap B(\widetilde{x}_0,D)=\{p\}$.
It follows that the exterior unit normal $\widetilde{w}_0\in T_p$
of $B(\widetilde{x}_0,D)$ at $p$ satisfies 
$\angle(u,o,\widetilde{w}_0)\leq \frac{\pi}2$, contradicting \eqref{uNtildeN}, and verifying \eqref{NtildeN}.

Finally to prove that $K$ is of constant width $D$, we choose any unit normal $\nu\in T_p$ to $K$ at $p$. We deduce from (iii) and \eqref{NtildeN} that $\nu$ is also normal to a $B(x,D)$ with $x\in K$ and $p\in \partial B(x,D)$, completing the proof of Proposition~\ref{completeness}.
\hfill $\Box$\\

The following Lemma is a corollary of the fact that bodies of constant width are exactly the complete sets in $\mathcal{M}^n$ (see of Proposition~\ref{completeness}).

\begin{lemma}\label{ConstantDiameter}
Let $\mathcal{M}^n$ be either $\R^n$, $S^n$ or $H^n$, and let
$K\subset\mathcal{M}^n$ be a complete set of diameter $D>0$ where $D<\frac{\pi}{2}$ if $\mathcal{M}^n=S^n$.
If  $B_{\mathcal{M}^n}\left(z,r\right)\subset K$ for $r>0$
with $y\in\partial_{\mathcal{M}^n}B_{\mathcal{M}^n}\left(z,r\right)\cap\partial_{\mathcal{M}^n}K$, then $w\in \partial K$ where $w$ is the unique point on the geodesic line through $y$ and $z$ such that $z\in\left[w,y\right]_{\mathcal{M}^n}$ and $d_{\mathcal{M}^n}\left(w,y\right)=D$.
\end{lemma}

Now we are ready to prove the following nice connection between the inradius and the circumradius of complete sets.

\begin{prop}\label{r_plus_R}
Let $\mathcal{M}^n$ be either $\R^n$, $S^n$ or $H^n$. If 
$K\subset\mathcal{M}^n$ is a complete set of diameter $D>0$ where $D<\frac{\pi}{2}$ in the case $\mathcal{M}^n=S^n$, then
$$
R\left(K\right)+r\left(K\right)=D.
$$
Furthermore, $K$ has a unique inscribed ball whose center is the circumcenter.
\end{prop}
\proof To simplify notation, we drop the reference to $\mathcal{M}^n$ in the formulas.

Let $p\in K$ be the unique circumcenter, and hence $K\subset B\left(p,R\left(K\right)\right)$. 
Since $R(K)<D$ by Lemma~\ref{jung}, and $K$ is of constant width $D$ by Proposition~\ref{completeness}, we deduce that $p\in{\rm int}\,K$.

Let $\varrho>0$ be the radius of the largest inscribed ball of $K$ with center $p$:
$$
\varrho=\max\left\{r>0\colon B\left(p,r\right)\subseteq K\right\}.
$$
 Let $z$ be a boundary point in $\partial K\cap\partial B\left(p,\varrho\right)$. By Lemma~\ref{ConstantDiameter} there is a $y\in\partial K$ such that $d\left(y,z\right)=D$ and $p\in\left[y,z\right]$, and hence
$$
R\left(K\right)\geq d\left(y,p\right)=D-\varrho,
$$
which in turn yields that
\begin{equation}\label{r_plus_R_lowerbound}
 r\left(K\right)\geq\varrho\geq D-R\left(K\right).   
\end{equation}
On the other hand, we claim if $B\left(w,r\right)\subset K$ for $r>0$ and $w\in K$, then
\begin{equation}\label{r_plus_R_upperbound}
r\leq D-R\left(K\right)
\end{equation}
with equality if and only if $w=p$ and $r=\varrho$.

According to Lemma~\ref{jung}, there exist
 some 
 $\left\{q_1;\ldots;q_k\right\}\subseteq\partial K\cap\partial B \left(p,R\left(K\right)\right)$ 
 for $2\leq k\leq n+1$ such that $p\in\left[q_1,\ldots,q_k\right]$ and $\left[q_1,\ldots,q_k\right]$ is a $\left(k-1\right)$-dimensional simplex in $\mathcal{M}^n$. 

If $w=p$, then 
$$
B\left(p,\varrho\right)\subset K\subset B\left(q_1,D\right),
$$
implying
$$
r\leq\varrho\leq D-d\left(p,q_1\right)=D-R\left(K\right).
$$
It follows from \eqref{r_plus_R_lowerbound} that we have $r=D-R\left(K\right)$ if and only if $r=\varrho$.

If $w\neq p$, then let $H$ be the hyperplane at $p$ orthogonal to $\left[w,p\right]$ and let $H^+$ be the closed half-space with $w\in{\rm{int}}H^+$. Since $p\not\in{\rm{int}}H^+$ and $p\in\left[q_1,\ldots,q_k\right]$, there is an index $1\leq i\leq k$ such that the point $q_i\not\in{\rm{int}}H^+$ by the convexity of the open half-space, and so $\angle\left(q_i,p,w\right)\geq\frac{\pi}{2}$. Now we have
$$
d\left(q_i,w\right)>d\left(q_i,p\right)=R\left(K\right)
$$
by the Law of Cosines for sides in $\mathcal{M}^n$.
It follows that
$$
B\left(w,r\right)\subset K\subset B\left(q_i,D\right),
$$
now implying
$$
r\leq D-d\left(q_i,w\right)<D-d\left(q_i,p\right)=
D-R\left(K\right),
$$
concluding the proof of \eqref{r_plus_R_upperbound}.

Finally, combining \eqref{r_plus_R_lowerbound}
and \eqref{r_plus_R_upperbound}
yields Proposition~\ref{r_plus_R}.
\hfill $\Box$\\

Combining Lemma~\ref{jung} and Proposition~\ref{r_plus_R} implies the following.

\begin{corollary}
\label{CompleteRadii}
Let $\mathcal{M}^n$ be either $\R^n$, $S^n$ or $H^n$. If 
$K\subset\mathcal{M}^n$ is a complete set of diameter $D>0$ where $D<\frac{\pi}{2}$ in the case $\mathcal{M}^n=S^n$, then
$$
R\left(K\right)\leq R\left(\Delta_{\mathcal{M}^n}(D)\right) \mbox{ \ and \ }r\left(K\right)\geq D-R\left(\Delta_{\mathcal{M}^n}(D)\right).
$$
Furthermore, either $R\left(K\right)= R\left(\Delta_{\mathcal{M}^n}(D)\right)$  or $r\left(K\right)= D-R\left(\Delta_{\mathcal{M}^n}(D)\right)$
if and only if $K$ contains a congruent copy of
$\Delta_{\mathcal{M}^n}(D)$.
\end{corollary}
\noindent{\bf Remark} According to Lemma~\ref{Zorn}, there exists a complete set $K_0$ of diameter $D$ in $\mathcal{M}^n$ containing a congruent copy of
$\Delta_{\mathcal{M}^n}(D)$. If $n=2$, then there exists a unique such set up to congruency; namely, the Reuleaux triangle (see Section~\ref{secReuleaux}).

\section{Complete sets in the hyperbolic space}
\label{sechyp}

If $\mathcal{M}^n=H^n$, we say that $K$ is of \emph{constant horospherical width} $D$ if $w_i\left(K\right)=D$ for every $i$ ideal point. The following lemma characterizes the completeness of a convex body in all spaces of constant curvature.

In the hyperbolic space $H^n$, we also consider horoballs and horospheres ``centered'' at ideal points. An ideal point of $H^n$ is represented by the linear hull of a $w\in \R^{n+1}$ with $\mathcal{B}_{H^n}(w,w)=0$ and  
$\mathcal{B}_{H^n}(w,e)>0$. Any hyperbolic line $\ell=\Pi\cap H^n$ for a linear two-plane $\Pi$ with $\Pi\cap H^n\neq\emptyset$ contains exactly two ideal points. If $p\in\ell$ and $v\in T_p\cap \Pi$ is one of the two tangent vectors to $\ell$ at $p$ with $\mathcal{B}_{H^n}(v,v)=-1$, then the two ideal points are represented by the spans of $u-v$ (``point of $\ell$ at infinity'' in the direction of $v$) and $u+v$ (``point of $\ell$ at infinity'' in the direction of $-v$).
 We also observe that hyperbolic lines containing either ideal point of $\ell$ are the lines parallel to $\ell$.

Let $w\in \R^{n+1}$ with $\mathcal{B}_{H^n}(w,w)=0$ and  
$\mathcal{B}_{H^n}(w,e)>0$, and let $i$ be the ideal point represented by ${\rm lin}\,w$. For any $p\in H^n$, there exists a unique line $\ell$ passing through $p$ and $i$, and a unit tangent vector $v\in T_p$ to $\ell$ is 
$u-\mathcal{B}_{H^n}(p,w)^{-1}w$.

For $s>0$, a horoball at $i$ is
$$
A=\{z\in H^n\colon\mathcal{B}_{H^n}(z,w)\geq s\},
$$
and the corresponding horosphere is
$$
\partial A=\{z\in H^n\colon\mathcal{B}_{H^n}(z,w)=s\}.
$$
We observe that for any $p\in\partial A$, the line passing through $p$ and $i$ is orthogonal to the horosphere $\partial A$. If $n=2$, then a horosphere is also called a horocycle.

If $r>s$, then the horoball $C=\{z\in H^n\colon\mathcal{B}_{H^n}(z,w)\geq r\}\subset A$ satisfies that the distance of the parallel horospheres $\partial A$ and $\partial C$ is $d>0$ where
$\frac{r}{s}=\cosh d+\sinh d$. In particular,  any hyperbolic line $\ell$ containing $i$ satisfies that the distance between 
$\ell\cap \partial A$ and $\ell\cap \partial C$ is $d$, and we set $d=w\left(A\setminus {\rm int}C\right)$
be the width of closed region bounded by the horospheres $\partial A$ and $\partial C$.

For Lemma~\ref{ballinhoroball}, we also introduce the Poincar\'e disk model of the hyperbolic $n$-space. In the Poincar\'e disk model, the hyperbolic $H^n$ space is identified with the interior of the unit Euclidean ball $B^n$ in $\R^n$, and the set of ideal points are just $\partial B^n$. 
A hyperbolic line in the Poincar\'e disk model is the intersection of ${\rm int}B^n$ and a Euclidean circle that is orthogonal to $\partial B^n$ at the two intersection points.  
A horosphere at an ideal point $i\in \partial B^n$
is of the form $\partial G\setminus\left\{i\right\}$ for a Euclidean $n$-ball $G\subset B^n$ of radius less than one and touching $\partial B^n$ in $i$.
In addition hyperbolic $n$-balls, in the Poincar\'e model coincide with Euclidean $n$-balls contained in ${\rm int}\,B^n$. In particular, the following Lemma~\ref{ballinhoroball} is  a simple consequence of these properties of the Poincar\'e disk model.

\begin{lemma}
\label{ballinhoroball}
Any horoball $A$ in $H^n$ is closed and convex, and if $A$ and a ball $B_{H^n}(z,r)$, $r>0$, have a common exterior normal $u\in T_y$ for some $y\in \partial A\cap \partial B_{H^n}(z,r)$, then
$B_{H^n}(z,r)\subset A$ with $B_{H^n}(z,r)\cap A=\{y\}$.
\end{lemma}

Lemma~\ref{ballinhoroball} directly yields the following.

\begin{corollary}\label{horostripwidth}
If $A$ and $C$ are horoballs at an ideal point $i$ of $H^n$ with $C\subset A$, then  $d_{H^n}(x,y)\geq w\left(A\setminus {\rm int}C\right)$ for any $x\in\partial A$ and $y\in\partial C$ with equality if and only if the line $\ell$ passing through $x$ and $y$ has $i$ as an ideal point (and hence $\ell$ is orthogonal to  $\partial A$ and $\partial C$).
\end{corollary}

We deduce the following properties of complete sets from Proposition~\ref{completeness}  and Lemma~\ref{ballinhoroball}.

\begin{corollary}\label{CompleteProperties}
Let $K\subset\mathcal{M}^n$ a $D$-complete set where $D>0$ and $D<\frac{\pi}{2}$ in the case of $\mathcal{M}^n=S^n$. Then the following hold:
\begin{enumerate}
    \item $K$ is convex and compact
    \item For all $z\in\partial_{\mathcal{M}^n}K$ there is a $y\in\partial_{\mathcal{M}^n}K$ such that $d_{\mathcal{M}^n}\left(y,z\right)=D$
    \item For for any pair of points $y,z\in K$ if $\sigma$ denotes the shorter arc of some circle of radius at least $D$ through $y$ and $z$, then $\sigma\subseteq K$
\end{enumerate}
\end{corollary}

For a bounded set $X\subset H^n$ and an ideal point $i$, we define the horospherical width corresponding to $i$ as
$$
w_i\left(X\right)=\min\left\{w\left(A\setminus {\rm int}A_0\right)\colon A\text{ and }A_0\text{ are horoballs at }i\text{ such that }X\subset A\setminus A_0\right\}.
$$
In the hyperbolic case, completeness can be also characterized in terms of constant horospherical width
analogously to the characterization of a complete set $Z$ in the Euclidean space in terms of the constant width of the parallel strips containing $Z$.

\begin{prop}\label{CompletenessHyp}
For $D>0$ and a compact convex set $K\subseteq H^n$, $K$ is constant width $D$ if and only if  $w_i\left(K\right)=D$ for any ideal point $i$.
\end{prop}
\proof 
First we assume that $K$ is of constant width $D$, and let $i$ be an ideal point. We consider the minimal horoball $A$ at $i$ containing $K$, and  
 the maximal horoball $A_0$ at $i$ with $K\cap{\rm int}A_0=\emptyset$. Let $x\in\partial K\cap \partial A$ and $y\in\partial K\cap \partial A_0$.
 It follows that $A$ and $K$ have a common exterior normal at $x$, and as $K$ is of constant width $D$,
 the line $\ell$ orthogonal to $A$ at $x$ (and hence containing $i$) intersects $\partial K$ in a point $z$ such that $d(x,z)=D$. We deduce from the definition of
 $w(A\cap {\rm int}A_0)$ and Corollary~\ref{horostripwidth} that
 $$
 D=d(x,z)\leq w(A\cap {\rm int}A_0)\leq d(x,y)\leq{\rm diam}\,K=D; 
 $$
therefore, $w(A\cap {\rm int}A_0)=D$.

Next we assume that $w_i\left(K\right)=D$ for any ideal point $i$. We consider $p,q\in \partial K$ such that
${\rm diam}\,K=d(p,q)$. Let $\ell$ be the line passing through $p$ and $q$, and let $i$ the ideal point
of $\ell$ such that $p$ lies between $i$ and $q$. Since
$K\subset B(p,D)\cap B(q,D)$, it follows from Lemma~\ref{ballinhoroball} that
${\rm diam}\,K=w(A\cap {\rm int}A_0)=D$
for the horoballs $A$ and $A_0$ at $i$ with
$p\in\partial A_0$ and $q\in\partial A$.

 Let $x\in\partial K$, let $u\in T_x$ be an exterior normal to $K$ at $x$. We consider the line $\ell$ passing through $x$ along $u$, and let $i$ the ideal point of $\ell$ in the direction of $u$, and hence the open half-line connecting $x$ to $i$ does not intersect $K$.  It follows from Lemma~\ref{ballinhoroball} that the tangent hyperplane $H$ at $x$ with normal $u$ separates $K$ and the horoball $\Omega_0$ at $i$ with $x\in\partial \Omega_0$. Let $\Omega$ be the minimal horoball at $i$ containing $K$, and let $y\in \partial\Omega\cap \partial K$. We deduce from  Corollary~\ref{horostripwidth} that
 $$
 D={\rm diam}\,K\geq d(x,y)\geq 
 w(\Omega\cap {\rm int}\Omega_0)= D; 
 $$
therefore, $d(x,y)=D$, and hence the equality case of 
Corollary~\ref{horostripwidth} yields that $y\in \ell$.
\hfill $\Box$\\

\noindent{\bf Proof of Theorem~\ref{CompletenessHyptheo} } 
Combining Proposition~\ref{completeness} and Proposition~\ref{CompletenessHyp} yields 
Theorem~\ref{CompletenessHyptheo}. \hfill$\Box$\\

\section{About the Reuleaux triangle} 
\label{secReuleaux}

Let $\mathcal{M}^2$ be either $\R^2$, $S^2$ or $H^2$.
If $v_1,v_2,v_3$ are the vertices of $\Delta_{\mathcal{M}^2}\left(D\right)$, then
a Reuleaux triangle $K$ of diameter $D$ is congruent to
$$
U_D=\bigcap_{i=1}^3{B_{\mathcal{M}^2}\left(v_i,D\right)}.
$$
We observe that $U_D$ is the unique complete set of diameter $D$ containing $\Delta_{\mathcal{M}^2}\left(D\right)$ according to Proposition~\ref{completeness}.
The following statement bounding the inradii and circumradii of two-dimensional complete sets is a direct consequence of
Corollary~\ref{CompleteRadii}.

\begin{prop}\label{reuleauxradii}
If $K\subseteq \mathcal{M}^2$ is complete with diameter $D>0$, then
$$
r\left(U_D\right)\leq r\left(K\right)\leq R\left(K\right)\leq R\left(U_D\right),
$$
and $r\left(U_D\right)=r\left(K\right)$ or $R\left(K\right)=R\left(U_D\right)$ yields that $K$ is a Reuleaux triangle.
\end{prop}

When discussing convex domains bounded by circular arcs, we need the following two elementary properties of circles.

\begin{lemma}
\label{DCircles}
Let $\mathcal{M}^2$ be either $H^2$, $\R^2$ or $S^2$.
For $D>0$ (where $D<\frac{\pi}2$ provided $\mathcal{M}^2=S^2$),
\begin{description}
\item[(i)] if $0<\varrho\leq D$ and $0<d_{\mathcal{M}^2}(x,y)<D+\varrho$, then 
$\partial B_{\mathcal{M}^2}(x,D)$ and 
$\partial B_{\mathcal{M}^2}(y,\varrho)$ intersect in two points $a$ and $b$, and if even $\varrho= D$, then the open longer
arc $\arc{ab}$ of 
$\partial B_{\mathcal{M}^2}(x,D)$ lies outside of 
$B_{\mathcal{M}^2}(y,D)$;
\item[(ii)] if $0<\varrho< D$ and $\varrho<d_{\mathcal{M}^2}(p,q)<2D-\varrho$, then
there exist exactly two circular discs
 of radius $D$ that contain $B_{\mathcal{M}^2}(p,\varrho)$
 and whose boundary touch $B_{\mathcal{M}^2}(p,\varrho)$ and contain $q$.
 \end{description}
\end{lemma}

Let $\mathcal{M}^2$ be either $H^2$, $\R^2$ or $S^2$.
We say that $\Gamma\subseteq \mathcal{M}^2$ is
 Jordan domain if is the union $\Gamma=\gamma\cup\Gamma_0$ where $\gamma\subseteq \mathcal{M}^2$ is a simple closed curve and $\Gamma_0$ is the bounded component 
 (or a bounded component if $\mathcal{M}^2=S^2$) of $\mathcal{M}^2\setminus\gamma$ provided by the Jordan curve theorem. In particular ${\rm{int}}\Gamma=\Gamma_0$ and $\partial_{\mathcal{M}^2}\Gamma=\gamma$.
 In order to specify the component of
 $S^2\setminus\gamma$ in the case
  $\mathcal{M}^2=S^2$, if $\gamma$ is contained in an open hemisphere $S_+$, then we assume that $\Gamma_0\subset S_+$.

 For $D>0$ where $D<\frac{\pi}2$ if $\mathcal{M}^2=S^2$, let 
$\Delta_{\mathcal{M}^2}\left(D\right)=\left[v_1,v_2,v_3\right]_{\mathcal{M}^2}$ be a regular triangle in $\mathcal{M}^2$. We write $\widetilde{p}$ to denote common circumcenter of $\Delta_{\mathcal{M}^2}\left(D\right)$ and the Reuleaux triangle 
$U_D=\cap_{i=1,2,3}B_{\mathcal{M}^2}(v_i,D)$ associated to $\Delta_{\mathcal{M}^2}\left(D\right)$. 

If $r\left(U_D\right)\leq\varrho\leq\frac{D}{2}$, then we define the compact convex sets $\widetilde{Q}_D\left(\varrho\right)$ and
$\widetilde{W}_D\left(\varrho\right)$ as follows. 
For $i=1,2,3$, let 
$\widetilde{q}_i=\widetilde{q}_i(\varrho,D)\in[\widetilde{p},v_i]_{\mathcal{M}^2}$ satisfy $d_{\mathcal{M}^2}(\widetilde{p},\widetilde{q}_i)=D-\varrho$,
and let $\widetilde{t}_i=\widetilde{t}_i(\varrho,D)$ be the point of the line passing through
$v_i,\widetilde{q}_i,\widetilde{p}$ such that
$\widetilde{p}\in[\widetilde{t}_i,\widetilde{q}_i]_{\mathcal{M}^2}$
and $d_{\mathcal{M}^2}(\widetilde{p},\widetilde{t}_i)=\varrho$.

To define $\widetilde{W}_D\left(\varrho\right)$,
we consider the circular arc $\sigma_{ij}$ of radius $D$ connecting $\widetilde{q}_i$ and $\widetilde{t}_j$
for $i\neq j$, $i,j\in\{1,2,3\}$ where the line through
$\widetilde{q}_i$ and $\widetilde{t}_j$ separates $\sigma_{ij}$ from
$\widetilde{p}$, and let
 $\widetilde{W}_D\left(\varrho\right)$ be the convex compact set  bounded by the six circular arcs $\sigma_{ij}$, $i\neq j$. In particular, $\widetilde{W}_D\left(r_K(U_D)\right)=U_D$.
 
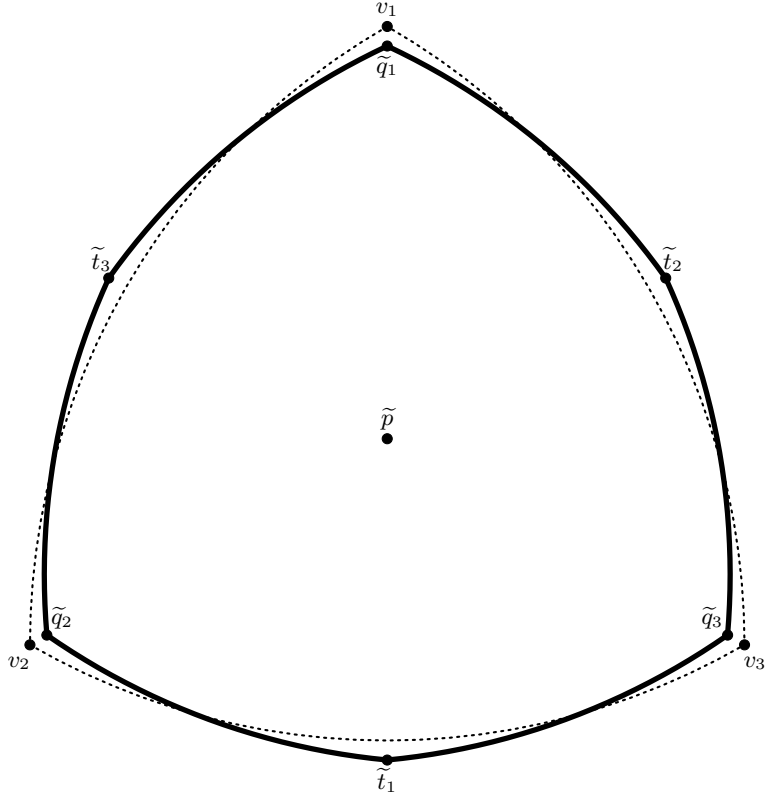
\begin{center}
\begin{tikzpicture}[line cap=round,line join=round,>=triangle 45,x=1cm,y=1cm, scale=9.5]
\clip(-0.072,-0.22) rectangle (1.144,0.9208701508801388);
\draw [shift={(1,0)},line width=.8pt,dotted]  plot[domain=2.0943951023931953:3.141592653589793,variable=\t]({1*1*cos(\t r)+0*1*sin(\t r)},{0*1*cos(\t r)+1*1*sin(\t r)});
\draw [shift={(0,0)},line width=.8pt,dotted]  plot[domain=0:1.0471975511965976,variable=\t]({1*1*cos(\t r)+0*1*sin(\t r)},{0*1*cos(\t r)+1*1*sin(\t r)});
\draw [shift={(0.5,0.8660254037844388)},line width=.8pt,dotted]  plot[domain=4.1887902047863905:5.235987755982988,variable=\t]({1*1*cos(\t r)+0*1*sin(\t r)},{0*1*cos(\t r)+1*1*sin(\t r)});
\draw [shift={(0.9246504784333377,-0.06668223747786052)},line width=2pt]  plot[domain=2.0093721182089483:2.5224251945822695,variable=\t]({1*1*cos(\t r)+0*1*sin(\t r)},{0*1*cos(\t r)+1*1*sin(\t r)});
\draw [shift={(0.5954232724203452,0.834111922683636)},line width=2pt]  plot[domain=4.103767220602143:4.616820296975464,variable=\t]({1*1*cos(\t r)+0*1*sin(\t r)},{0*1*cos(\t r)+1*1*sin(\t r)});
\draw [shift={(1.020073750853683,0.09859571857866334)},line width=2pt]  plot[domain=2.7135625614007193:3.2266156377740405,variable=\t]({1*1*cos(\t r)+0*1*sin(\t r)},{0*1*cos(\t r)+1*1*sin(\t r)});
\draw [shift={(0.4045767275796549,0.8341119226836362)},line width=2pt]  plot[domain=4.807957663793914:5.321010740167235,variable=\t]({1*1*cos(\t r)+0*1*sin(\t r)},{0*1*cos(\t r)+1*1*sin(\t r)});
\draw [shift={(-0.02007375085368288,0.09859571857866313)},line width=2pt]  plot[domain=-0.08502298418424736:0.4280300921890742,variable=\t]({1*1*cos(\t r)+0*1*sin(\t r)},{0*1*cos(\t r)+1*1*sin(\t r)});
\draw [shift={(0.07534952156666211,-0.06668223747786037)},line width=2pt]  plot[domain=0.6191674590075236:1.1322205353808448,variable=\t]({1*1*cos(\t r)+0*1*sin(\t r)},{0*1*cos(\t r)+1*1*sin(\t r)});
\begin{scriptsize}
\draw [fill=black] (0,0) circle (.2pt);
\draw[color=black] (-0.015,-0.025) node {$v_2$};
\draw [fill=black] (1,0) circle (.2pt);
\draw[color=black] (1.015,-0.025) node {$v_3$};
\draw [fill=black] (0.5,0.8660254037844388) circle (.2pt);
\draw[color=black] (0.5,0.89) node {$v_1$};
\draw [fill=black] (0.5,0.2886751345948129) circle (.2pt);
\draw[color=black] (0.5,0.31614332483983393) node {$\widetilde{p}$};
\draw [fill=black] (0.5,0.838675134594813) circle (.2pt);
\draw[color=black] (0.5,0.81) node {$\widetilde{q}_1$};
\draw [fill=black] (0.5,-0.1613248654051871) circle (.2pt);
\draw[color=black] (0.5,-0.189) node {$\widetilde{t}_1$};
\draw [fill=black] (0.9763139720814413,0.013675134594812843) circle (.2pt);
\draw[color=black] (0.955,0.038) node {$\widetilde{q}_3$};
\draw [fill=black] (0.11028856829700262,0.513675134594813) circle (.2pt);
\draw[color=black] (0.1,0.54) node {$\widetilde{t}_3$};
\draw [fill=black] (0.023686027918558724,0.013675134594812873) circle (.2pt);
\draw[color=black] (0.045,0.038) node {$\widetilde{q}_2$};
\draw [fill=black] (0.8897114317029974,0.513675134594813) circle (.2pt);
\draw[color=black] (0.9,0.54) node {$\widetilde{t}_2$};
\end{scriptsize}
\end{tikzpicture}
\end{center}
 
For $\widetilde{Q}_D\left(\varrho\right)$, we set
$\widetilde{Q}_D\left(r(U_D)\right)=U_D$ when $\varrho=r(U_D)$ and $\widetilde{Q}_D\left(D/2\right)=U_D$ when $\varrho=D/2$. Next let $r\left(U_D\right)<\varrho<\frac{D}{2}$.
For each $\widetilde{q}_i$, $i=1,2,3$, 
Lemma~\ref{DCircles} (ii)
provides two circular disks 
$B_{\mathcal{M}^2}(\widetilde{z}_{ij},D)$ and
$B_{\mathcal{M}^2}(\widetilde{z}_{ik},D)$ containing
$B_{\mathcal{M}^2}(\widetilde{p},\varrho)$
with $\left\{i,j,k\right\}=\left\{1,2,3\right\}$ such that $\widetilde{q}_i\in\partial B_{\mathcal{M}^2}\left(\widetilde{z}_{ij},D\right)$ and $B_{\mathcal{M}^2}\left(\widetilde{p},\varrho\right)$ touches 
$\partial B_{\mathcal{M}^2}\left(\widetilde{z}_{ij},D\right)$
in $\widetilde{w}_{ij}$ where $\widetilde{w}_{ij}$ and $\widetilde{q}_j$ lie on the same side of the geodesic line through $\widetilde{q}_i$ and $\widetilde{t}_i$. 

Next, if $\left\{i,j,k\right\}=\left\{1,2,3\right\}$,
then $\widetilde{q}_j,\widetilde{q}_k\in{\rm int}B_{\mathcal{M}^2}\left(\widetilde{q}_{i},D\right)$
and $\partial B_{\mathcal{M}^2}\left(\widetilde{q}_{i},D\right)$
touches $B_{\mathcal{M}^2}\left(\widetilde{p},\varrho\right)$ in $\widetilde{t}_i$. Therefore, $\widetilde{t}_i$ is contained
in the open shorter arc $\arc{\widetilde{w}_{jk}\widetilde{w}_{kj}}$ of 
$\partial B_{\mathcal{M}^2}\left(\widetilde{p},\varrho\right)$.
In particular, we can define $\widetilde{Q}_D\left(\varrho\right)$
as the compact convex set in $\mathcal{M}^2$ that is bounded
by the six shorter circular arcs of the from 
$\arc{\widetilde{q}_{i}\widetilde{w}_{ij}}$ of
$\partial B_{\mathcal{M}^2}\left(\widetilde{z}_{ij},D\right)$
for $i,j=1,2,3$, $i\neq j$, and the three 
shorter circular arcs of the from 
$\arc{\widetilde{w}_{ij}\widetilde{w}_{ji}}$
of $B_{\mathcal{M}^2}\left(\widetilde{p},\varrho\right)$
for $1\leq i<j\leq 3$.

For $\{i,j,k\}=\{1,2,3\}$, let $\widetilde{\Gamma}_i\left(\varrho\right)$ denote the Jordan domain enclosed by the shorter arc  $\arc{\widetilde{w}_{ij}\widetilde{q}_i}$ of  $\partial B_{\mathcal{M}^2}\left(\widetilde{z}_{ij},D\right)$; the shorter arc $\arc{\widetilde{w}_{ik}\widetilde{q}_i}$ of  $\partial B_{\mathcal{M}^2}\left(\widetilde{z}_{ik},D\right)$, and the shorter arc $\arc{\widetilde{w}_{ij}\widetilde{w}_{ik}}$ of $\partial B_{\mathcal{M}^2}\left(\widetilde{p},\varrho\right)$. Then
\begin{eqnarray}
\nonumber
\widetilde{Q}_D\left(\varrho\right)&=&
B_{\mathcal{M}^2}\left(\widetilde{p},\varrho\right)\cup\widetilde{\Gamma}_1\left(\varrho\right)\cup\widetilde{\Gamma}_2\left(\varrho\right)\cup\widetilde{\Gamma}_3\left(\varrho\right),\\
\label{VQ}
V\left(\widetilde{Q}_D\left(\varrho\right)\right)&=&
V\left(B_{\mathcal{M}^2}\left(\widetilde{p},\varrho\right)\right)+
V\left(\widetilde{\Gamma}_1\left(\varrho\right)\right)+
V\left(\widetilde{\Gamma}_2\left(\varrho\right)\right)+
V\left(\widetilde{\Gamma}_3\left(\varrho\right)\right)
\end{eqnarray}
where $\widetilde{\Gamma}_1\left(\varrho\right)$, $\widetilde{\Gamma}_2\left(\varrho\right)$ and $\widetilde{\Gamma}_3\left(\varrho\right)$ are congruent.

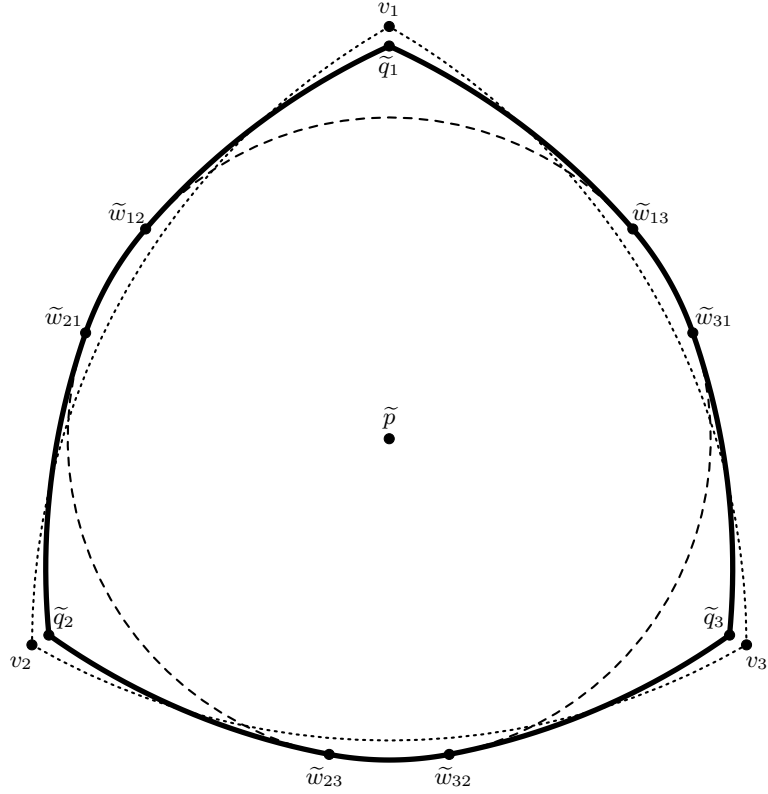
\begin{center}
\begin{tikzpicture}[line cap=round,line join=round,>=triangle 45,x=1cm,y=1cm, scale=9.5]
\clip(-0.072,-0.22) rectangle (1.144,0.9208701508801388);
\draw [shift={(1,0)},line width=.8pt,dotted]  plot[domain=2.0943951023931953:3.141592653589793,variable=\t]({1*1*cos(\t r)+0*1*sin(\t r)},{0*1*cos(\t r)+1*1*sin(\t r)});
\draw [shift={(0,0)},line width=.8pt,dotted]  plot[domain=0:1.0471975511965976,variable=\t]({1*1*cos(\t r)+0*1*sin(\t r)},{0*1*cos(\t r)+1*1*sin(\t r)});
\draw [shift={(0.5,0.8660254037844388)},line width=.8pt,dotted]  plot[domain=4.1887902047863905:5.235987755982988,variable=\t]({1*1*cos(\t r)+0*1*sin(\t r)},{0*1*cos(\t r)+1*1*sin(\t r)});
\draw [shift={(0.3973170456844901,0.8290048588308937)},line width=2pt]  plot[domain=4.900187199278438:5.329886865429862,variable=\t]({1*1*cos(\t r)+0*1*sin(\t r)},{0*1*cos(\t r)+1*1*sin(\t r)});
\draw [shift={(-0.019280744766041115,0.10743631944964122)},line width=2pt]  plot[domain=-0.09389910944687418:0.3358005567045506,variable=\t]({1*1*cos(\t r)+0*1*sin(\t r)},{0*1*cos(\t r)+1*1*sin(\t r)});
\draw [shift={(0.916597790450531,-0.0704157744960961)},line width=2pt]  plot[domain=2.0004959929463215:2.430195659097746,variable=\t]({1*1*cos(\t r)+0*1*sin(\t r)},{0*1*cos(\t r)+1*1*sin(\t r)});
\draw [shift={(0.08340220954946884,-0.07041577449609593)},line width=2pt]  plot[domain=0.7113969944920471:1.1410966606434718,variable=\t]({1*1*cos(\t r)+0*1*sin(\t r)},{0*1*cos(\t r)+1*1*sin(\t r)});
\draw [shift={(0.6026829543155101,0.8290048588308936)},line width=2pt]  plot[domain=4.094891095339516:4.524590761490941,variable=\t]({1*1*cos(\t r)+0*1*sin(\t r)},{0*1*cos(\t r)+1*1*sin(\t r)});
\draw [shift={(1.0192807447660412,0.10743631944964124)},line width=2pt]  plot[domain=2.805792096885243:3.2354917630366673,variable=\t]({1*1*cos(\t r)+0*1*sin(\t r)},{0*1*cos(\t r)+1*1*sin(\t r)});
\draw [shift={(0.5,0.2886751345948129)},line width=2pt]  plot[domain=2.4301956590977465:2.805792096885243,variable=\t]({1*0.45*cos(\t r)+0*0.45*sin(\t r)},{0*0.45*cos(\t r)+1*0.45*sin(\t r)});
\draw [shift={(0.5,0.2886751345948129)},line width=.8pt,dashed]  plot[domain=2.805792096885243:4.524590761490941,variable=\t]({1*0.45*cos(\t r)+0*0.45*sin(\t r)},{0*0.45*cos(\t r)+1*0.45*sin(\t r)});
\draw [shift={(0.5,0.2886751345948129)},line width=2pt]  plot[domain=4.524590761490941:4.900187199278438,variable=\t]({1*0.45*cos(\t r)+0*0.45*sin(\t r)},{0*0.45*cos(\t r)+1*0.45*sin(\t r)});
\draw [shift={(0.5,0.2886751345948129)},line width=.8pt,dashed]  plot[domain=-1.3829981079011482:0.33580055670455067,variable=\t]({1*0.45*cos(\t r)+0*0.45*sin(\t r)},{0*0.45*cos(\t r)+1*0.45*sin(\t r)});
\draw [shift={(0.5,0.2886751345948129)},line width=2pt]  plot[domain=0.3358005567045506:0.7113969944920473,variable=\t]({1*0.45*cos(\t r)+0*0.45*sin(\t r)},{0*0.45*cos(\t r)+1*0.45*sin(\t r)});
\draw [shift={(0.5,0.2886751345948129)},line width=.8pt,dashed]  plot[domain=0.7113969944920473:2.4301956590977465,variable=\t]({1*0.45*cos(\t r)+0*0.45*sin(\t r)},{0*0.45*cos(\t r)+1*0.45*sin(\t r)});
\begin{scriptsize}
\draw [fill=black] (0,0) circle (.2pt);
\draw[color=black] (-0.015,-0.025) node {$v_2$};
\draw [fill=black] (1,0) circle (.2pt);
\draw[color=black] (1.015,-0.025) node {$v_3$};
\draw [fill=black] (0.5,0.8660254037844388) circle (.2pt);
\draw[color=black] (0.5,0.89) node {$v_1$};
\draw [fill=black] (0.5,0.2886751345948129) circle (.2pt);
\draw[color=black] (0.5,0.3185880999301591) node {$\widetilde{p}$};
\draw [fill=black] (0.5,0.838675134594813) circle (.2pt);
\draw[color=black] (0.5,0.81) node {$\widetilde{q}_1$};
\draw [fill=black] (0.9763139720814413,0.013675134594812843) circle (.2pt);
\draw[color=black] (0.955,0.038) node {$\widetilde{q}_3$};
\draw [fill=black] (0.023686027918558724,0.013675134594812873) circle (.2pt);
\draw[color=black] (0.045,0.038) node {$\widetilde{q}_2$};
\draw [fill=black] (0.5840133262581445,-0.15341282159834363) circle (.2pt);
\draw[color=black] (0.59,-0.184) node {$\widetilde{w}_{32}$};
\draw [fill=black] (0.9248660638994881,0.4369614378954079) circle (.2pt);
\draw[color=black] (0.955,0.46) node {$\widetilde{w}_{31}$};
\draw [fill=black] (0.15914726235865634,0.5824767874873746) circle (.2pt);
\draw[color=black] (0.1335,0.605) node {$\widetilde{w}_{12}$};
\draw [fill=black] (0.8408527376413436,0.5824767874873749) circle (.2pt);
\draw[color=black] (0.8665,0.605) node {$\widetilde{w}_{13}$};
\draw [fill=black] (0.4159866737418554,-0.15341282159834377) circle (.2pt);
\draw[color=black] (0.41,-0.184) node {$\widetilde{w}_{23}$};
\draw [fill=black] (0.07513393610051212,0.4369614378954078) circle (.2pt);
\draw[color=black] (0.045,0.46) node {$\widetilde{w}_{21}$};
\end{scriptsize}
\end{tikzpicture}
\end{center}

We deduce from Lemma~\ref{DCircles} (i) that
if $r\left(U_D\right)<\varrho\leq\frac{D}{2}$, then
\begin{equation}
\label{QDinWD}
\widetilde{W}_D\left(\varrho\right) \subsetneq
\widetilde{Q}_D\left(\varrho\right).
\end{equation}

\begin{lemma}
\label{diskdomains}
Using the notation as above in $\mathcal{M}^2$, if $r\left(U_D\right)<\varrho\leq\frac{D}{2}$, then
$$
V\left(\widetilde{W}_D\left(\varrho\right)\right)\geq
V\left(U_D\right).
$$
\end{lemma}
\proof To simplify notation, we drop the reference to
$\mathcal{M}^2$. For the proof, we keep  using the notation introduced for $\widetilde{W}_D\left(\varrho\right)$ as above.

For $i=1,2,3$, let $\{u_i\}=\partial B\left(\widetilde{p},r(U_D)\right)\cap \partial 
B\left(v_i,D\right)$ be the point such that
$\widetilde{p}\in[u_i,v_i]$, and for 
$i,j=1,2,3$, $i\neq j$, let $\gamma_{ij}\subset \partial U_D$ be the shorter arc $\arc{v_iu_j}$ of $\partial 
B\left(v_i,D\right)$. We frequently use the observation that
if $i\neq j$,
\begin{equation}
\label{ptpv}
d(\widetilde{p},\widetilde{t}_j)=\varrho<D-r(U_D)=
d(\widetilde{p},v_i).
\end{equation}

Since both $\widetilde{W}_D\left(\varrho\right)$ and $U_D$ are symmetric through the line $\ell_i$ passing through $v_i$, $\widetilde{q}_i$ $\widetilde{p}$, $u_i$ and $\widetilde{t}_i$ for
$i=1,2,3$, we concentrate on the part of
$\widetilde{W}_D\left(\varrho\right)$ and $U_D$ between
the half-lines $\widetilde{p}v_1$ and $\widetilde{p}\widetilde{t}_2$.
As 
$\widetilde{p}\in {\rm int}B(v_2,D)$, we deduce that
\begin{equation}
 \label{q1tildet2}
 \widetilde{q}_1\in {\rm int}B(v_2,D)\mbox{ \ and \ }
\widetilde{t}_2\not\in B(v_2,D).
\end{equation} 
Adding the fact
$\gamma_{12}\subset \partial B(v_2,D)$, it follows that $\gamma_{12}$ intersects $\sigma_{12}\subset \partial B(\widetilde{z}_{12},D)$ connecting
$\widetilde{q}_1$ and $\widetilde{t}_2$ in a point $d$ (and $\{d\}=\gamma_{12}\cap\sigma_{12}$  by Lemma~\ref{DCircles}).

Let $H_1$ be the closed half-plane bounded by $\ell_1$
and containing $u_2$ and $\widetilde{t}_2$, and let
$H_2$ be the closed half-plane bounded by $\ell_2$
and containing $v_1$ and $\widetilde{q}_1$. We consider the Jordan domains
\begin{eqnarray*}
\Omega_+&=&\mbox{closure of }\left(H_2\cap B(\widetilde{z}_{12},D)\right)\setminus B(v_2,D);\\
\Omega_-&=&\mbox{closure of }\left(H_1\cap B(v_2,D)\right)\setminus B(\widetilde{z}_{12},D).
\end{eqnarray*}
In particular, 
$\Omega_+$ is bounded by 
$[u_2,\widetilde{t}_2]$ and one subarc of
$\sigma_{12}$ and $\gamma_{12}$ each, and
$\Omega_-$ is bounded by 
$[v_1,\widetilde{q}_1]$ and one subarc of
$\sigma_{12}$ and $\gamma_{12}$ each, and 
$\Omega_+$ and $\Omega_-$ meet at $d$.
It follows from the symmetries 
$\widetilde{W}_D\left(\varrho\right)$ and $U_D$ 
\begin{equation}
\label{Omega+Omega-}
V\left(\widetilde{W}_D\left(\varrho\right)\right)-V(U_D)=
6\cdot\left(V\left(\Omega_+\right)-V(\Omega_-)\right).
\end{equation}

Let $\ell$ be the line perpendicularly bisecting the segment $[v_1,\widetilde{t}_2]$. It follows from \eqref{ptpv} that
the open half-plane bounded by $\ell$ and containing 
$\widetilde{t}_2$ also contains $\widetilde{p}$; therefore, the convexity of this half-plane yields that it also contains $u_2\in[\widetilde{p},\widetilde{t}_2]$. We deduce that $\ell$ intersects $\gamma_{12}$ in a point $e$. 
For any object $X$ in $\mathcal{M}^2$, we write $X'$ to denote the reflected image of $X$ through $\ell$. In particular, $\gamma'_{12}\subset \partial B(v'_2,D)$ 
is the shorter arc connecting $\widetilde{t}_2$ and $u'_2$.
It follows that $\gamma'_{12}\cap \sigma_{12}=\{e\}$, and 
$\gamma'_{12}$ is orthogonal to the line $\ell'_2$ containing $v_1$, $u'_2$, $\widetilde{p}'$ and $v'_2$ at $u'_2$. 

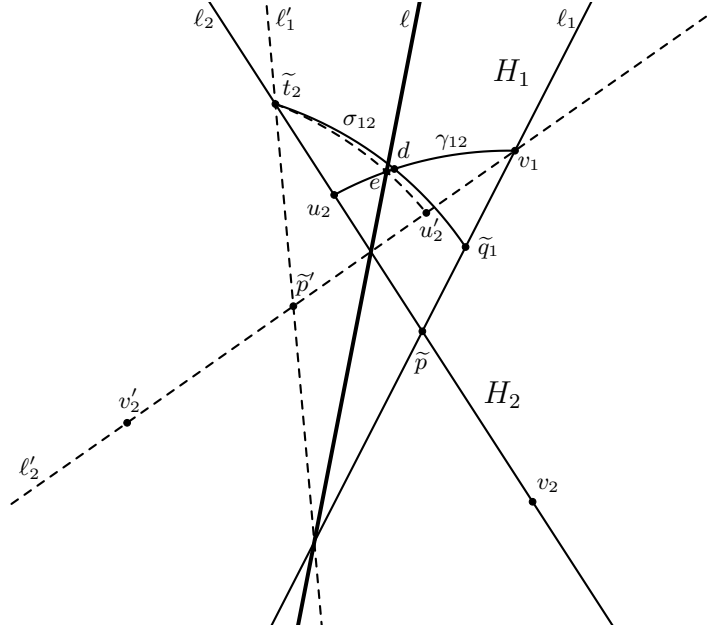
\begin{center}
\begin{tikzpicture}[line cap=round,line join=round,>=triangle 45,x=1cm,y=1cm, scale=1.8]
\clip(-5.6,-1.16) rectangle (-0.41810526494604466,3.4413791957373294);
\draw [line width=.8pt,domain=-5.9582408008551:-0.41810526494604466] plot(\x,{(--4.081649179682207--1.3349454680607475*\x)/0.6840472186216915});
\draw [line width=.8pt,domain=-5.9582408008551:-0.41810526494604466] plot(\x,{(-2.3778686974882723-1.259875002744846*\x)/0.8140730786966697});
\draw [line width=1.5pt,domain=-5.9582408008551:-0.41810526494604466] plot(\x,{(--5.71781084277351--1.7694779902172493*\x)/0.344887868932378});
\draw [line width=.8pt,domain=-5.9582408008551:-0.41810526494604466, dashed] plot(\x,{(-4.4784274555707775-0.8619274616690251*\x)/-1.2276322946309253});
\draw [shift={(-4.457997850935183,0.12140345506138658)},line width=.8pt]  plot[domain=0.5941940749558408:1.257329704572253,variable=\t]({1*2.7*cos(\t r)+0*2.7*sin(\t r)},{0*2.7*cos(\t r)+1*2.7*sin(\t r)});
\draw [shift={(-1.902465184544648,-0.3546538719065175)},line width=.8pt]  plot[domain=1.553568658646083:2.0684109867763683,variable=\t]({1*2.7*cos(\t r)+0*2.7*sin(\t r)},{0*2.7*cos(\t r)+1*2.7*sin(\t r)});
\draw [line width=.8pt,domain=-5.9582408008551:-0.41810526494604466, dashed] plot(\x,{(--5.060079024396513--1.4941244844910768*\x)/-0.1326349307093514});
\draw [shift={(-4.596162098180612,0.17037294008540232)},line width=.8pt, dashed]  plot[domain=1.553568658646083:2.0684109867763683,variable=\t]({-0.9268014509927892*2.7*cos(\t r)+0.3755516881038621*2.7*sin(\t r)},{0.3755516881038621*2.7*cos(\t r)+0.9268014509927892*2.7*sin(\t r)});
\draw (-2.1,3.1) node[anchor=north west] {$H_1$};
\draw (-2.170048263818811,0.7212571711716464) node[anchor=north west] {$H_2$};
\begin{scriptsize}
\draw [fill=black] (-2.54,1.01) circle (.7pt);
\draw[color=black] (-2.55,0.78) node {$\widetilde{p}$};
\draw [fill=black] (-1.7259269213033304,-0.24987500274484598) circle (.7pt);
\draw[color=black] (-1.61,-0.13) node {$v_2$};
\draw [fill=black] (-1.8559527813783085,2.3449454680607476) circle (.7pt);
\draw[color=black] (-1.75,2.25) node {$v_1$};
\draw[color=black] (-1.47,3.32) node {$\ell_{1}$};
\draw[color=black] (-4.15,3.32) node {$\ell_{2}$};
\draw [fill=black] (-3.1912584629573355,2.0179000021958764) circle (.7pt);
\draw[color=black] (-3.3,1.9) node {$u_{2}$};
\draw [fill=black] (-3.625430771595558,2.6898333369931255) circle (.7pt);
\draw[color=black] (-3.5,2.83) node {$\widetilde{t}_{2}$};
\draw[color=black] (-2.67,3.32) node {$\ell$};
\draw[color=black] (-5.44,0) node {$\ell_{2}'$};
\draw [fill=black] (-2.2207779646432106,1.6329745517616827) circle (.7pt);
\draw[color=black] (-2.05,1.63) node {$\widetilde{q}_{1}$};
\draw[color=black] (-2.32,2.41) node {$\gamma_{12}$};
\draw [fill=black] (-2.746830840386753,2.209921461210904) circle (.7pt);
\draw[color=black] (-2.68,2.35) node {$d$};
\draw [fill=black] (-2.8044592575428418,2.190224650199426) circle (.7pt);
\draw[color=black] (-2.89,2.1) node {$e$};
\draw [fill=black] (-3.4927958408862074,1.195708852502049) circle (.7pt);
\draw[color=black] (-3.399481947238297,1.3628678747062075) node {$\widetilde{p}'$};
\draw [fill=black] (-4.720428135517133,0.33378139083302394) circle (.7pt);
\draw[color=black] (-4.7,0.5) node {$v_{2}'$};
\draw[color=black] (-3.5454771971443604,3.32) node {$\ell_{1}'$};
\draw[color=black] (-3,2.56) node {$\sigma_{12}$};
\draw [fill=black] (-2.5106900051814676,1.885250821837269) circle (.7pt);
\draw[color=black] (-2.47,1.76) node {$u_{2}'$};
\end{scriptsize}
\end{tikzpicture}
\end{center}

We also consider the Jordan domain
$$
\widetilde{\Omega}=\mbox{closure of }\left(H_2\cap B(v'_2,D)\right)\setminus B(v_2,D),
$$
that is bounded by $[u_2,\widetilde{t}_2]$ and one subarc of
$\gamma'_{12}$ and $\gamma_{12}$ each meeting at $e$. It follows that
$$
\widetilde{\Omega}'=\mbox{closure of }\left(H'_2\cap B(v_2,D)\right)\setminus B(v'_2,D)
$$
is bounded by $[u'_2,v_1]$ and one subarc of
$\gamma'_{12}$ and $\gamma_{12}$ each meeting at $e$.

Next we claim that
\begin{equation}
\label{q1inOmega}
q_1\in{\rm int}\,\widetilde{\Omega}'.
\end{equation}
First, $q_1\in{\rm int}\,B(v_2,D)$ according to 
\eqref{q1tildet2}.
It follows from \eqref{ptpv} that
$\angle(\widetilde{p},v_1,\widetilde{t}_2)<\angle(\widetilde{p},\widetilde{t}_2,v_1)$, and hence
\begin{equation}
\label{q1angle}
\angle(q_1,v_1,\widetilde{t}_2)=\angle(\widetilde{p},v_1,\widetilde{t}_2)<\angle(\widetilde{p},\widetilde{t}_2,v_1)=
\angle(u_2,\widetilde{t}_2,v_1)=\angle(u'_2,v_1,\widetilde{t}_2).
\end{equation}
As $q_1\in{\rm int}\,B(v_2,D)$, we deduce that 
$q_1\in{\rm int}\,H'_2$.

Finally, $q_1\in \partial B\left(v_1,\varrho-r(U_D)\right)$
as $d(q_1,v_1)=d(\widetilde{t}_2,u_2)=\varrho-r(U_D)$.
On the other hand, $B\left(\widetilde{t}_2,\varrho-r(U_D)\right)$ touches $B(v_2,D)$ from the outside by construction; therefore, $B\left(v_1,\varrho-r(U_D)\right)=B\left(\widetilde{t}_2,\varrho-r(U_D)\right)'$ touches $B(v'_2,D)$ from the outside at $u'_2$. We conclude that $q_1\not\in B(v'_2,D)$, completing the proof of \eqref{q1inOmega}.

It follows from \eqref{q1inOmega} that $d$ is contained in the shorter arc $\arc{e,v_1}$ of
$\gamma_{12}\subset \partial B(V_2,D)$, which fact directly yields that
$\widetilde{\Omega}\subset \Omega_+$, and together
with \eqref{q1inOmega} implies that 
$\Omega_-\subset\widetilde{\Omega}'$ (where we have tacitly used Lemma~\ref{DCircles} (i) concerning intersection patterns of circular arcs). We deduce from \eqref{Omega+Omega-} that
$$
V\left(\widetilde{W}_D\left(\varrho\right)\right)-V(U_D)=
6\cdot\left(V\left(\Omega_+\right)-V(\Omega_-)\right)\geq
6\cdot\left(V(\widetilde{\Omega})-V(\widetilde{\Omega}')\right)=0,
$$
completing the proof of Lemma~\ref{diskdomains}.
\hfill $\Box$

\section{Proof of Theorem~\ref{BL-hyp}}
\label{sec-BL-hyp}

We are now ready to prove the Blaschke--Lebesgue--Leichtweiss theorem in any surface of constant curvature.

\begin{theorem}
\label{BL-M2} If $K\subset \mathcal{M}^2$ is a body of constant width $D$ for $D>0$ (where $D<\frac{\pi}2$ provided $\mathcal{M}^2=S^2$), and $U_D\subset \mathcal{M}^2$ is a 
Reuleaux triangle of width $D$, then
$$
V_{\mathcal{M}^2}\left(K\right)\geq V_{\mathcal{M}^2}\left(U_D\right)
$$
with equality if and only if $K$ is congruent with $U_D$.
\end{theorem}
\proof Since we are only using the intrinsic geometry of 
$\mathcal{M}^2$, we drop it from the notation.
Let $\varrho=r(K)$, and hence $R(K)=D-\varrho$ according to
Lemma~\ref{r_plus_R}. 
According to Proposition~\ref{reuleauxradii}, we may assume that  $\varrho> r(U_D)$.

Let $p$ be the unique circumcenter of $K$
(guaranteed by Lemma~\ref{jung}). Lemma~\ref{r_plus_R} states that $p$ is also the unique incenter of $K$, so we have
$$
B\left(p,\varrho\right)\subseteq K\subseteq B\left(p,D-\varrho\right).
$$
It follows from Lemma~\ref{closest-points-on-the-boundary} that there are points $t_1;\ldots;t_k\in\partial K\cap\partial B\left(p,\varrho\right)$ such that $\left[t_1,\ldots,t_k\right]$ is a $\left(k-1\right)$-dimensional simplex and $p\in{\rm{relint}}\left[t_1,\ldots,t_k\right]$, where $k\in\left\{2;3\right\}$. In addition, there are points $q_1;\ldots;q_k\in\partial K$ such that $d\left(t_j,q_j\right)=D$ and $p\in\left[t_j,q_j\right]$ for all $1\leq j\leq k$. The rest of the proof will be divided into two cases depending on $k$.\\

\noindent{\bf Case 1 }$k=2$

If there are boundary points $t_1;t_2\in\partial K$ such that $p\in\left[t_1,t_2\right]$, then $q_1=t_2$ and $q_2=t_1$, so $\varrho=\frac{D}{2}$, and hence $K=B\left(p,\frac{D}{2}\right)$. It follows from the  isodiametric inequality in $\mathcal{M}^2$ (see B\"or\"oczky, Sagmeister \cite{BoS20} or Schmidt \cite{Sch48,Sch49}) that
$$
V\left(B\left(p,\frac{D}{2}\right)\right)>
V\left(U_D\right).
$$

\noindent{\bf Case 2 }$k=3$

Let $\left\{i,j,k\right\}=\left\{1,2,3\right\}$ and let $B\left(z_{ij},D\right)$ be the disk such that $B\left(p,\varrho\right)\subseteq B\left(z_{ij},D\right)$, $q_i\in\partial B\left(z_{ij},D\right)$ and $B\left(z_{ij},D\right)$ touches $B\left(p,\varrho\right)$ in $w_{ij}$ where 
$w_{ij}$ and $q_j$ lie on the same side of
the geodesic line through $q_i$ and $t_i$. In addition, Lemma~\ref{closest-points-on-the-boundary} yields that \begin{equation}
\label{pintt1t2t3}
p\in{\rm{int}}\left[t_1,t_2,t_3\right]
\mbox{ \ and \ }
p\in{\rm{int}}\left[q_1,q_2,q_3\right].
\end{equation}
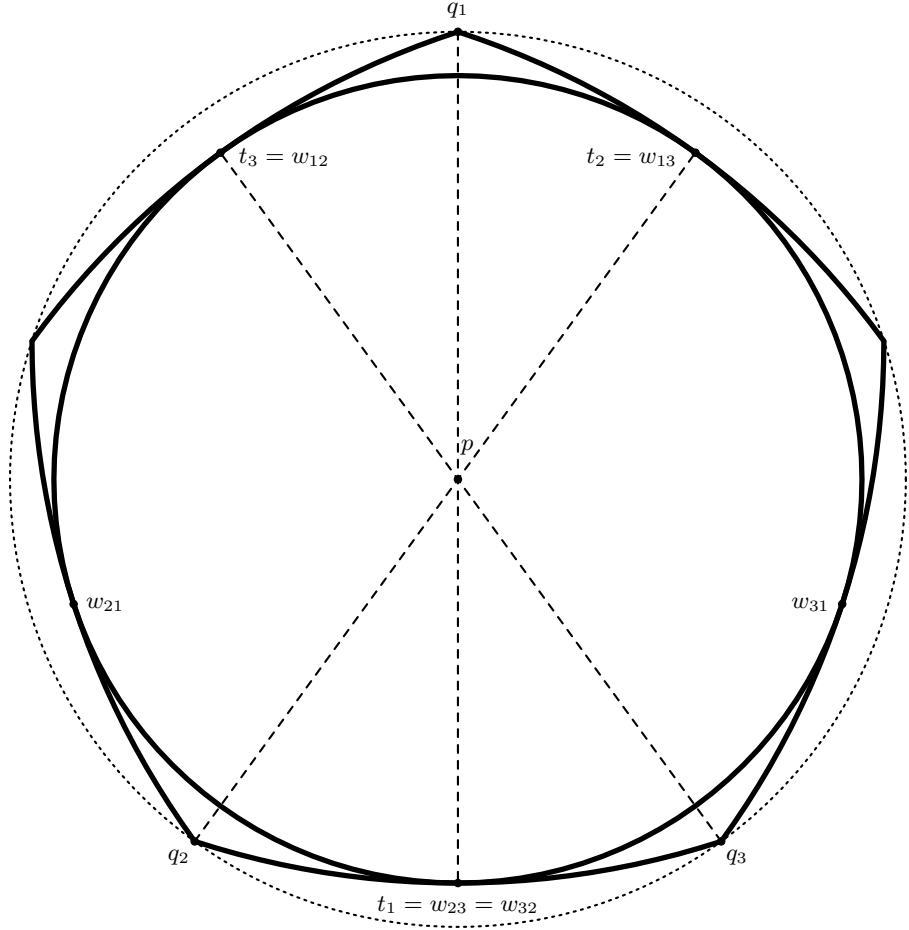
\begin{center}
\begin{tikzpicture}[line cap=round,line join=round,>=triangle 45,x=1cm,y=1cm, scale=7]
\clip(-0.5,-0.5) rectangle (1.5,1.8);
\draw [shift={(0.5,1.5388417685876266)},line width=2pt]  plot[domain=4.39822971502571:5.026548245743669,variable=\t]({1*1.6180339887498947*cos(\t r)+0*1.6180339887498947*sin(\t r)},{0*1.6180339887498947*cos(\t r)+1*1.6180339887498947*sin(\t r)});
\draw [shift={(-0.30901699437494734,0.9510565162951536)},line width=2pt]  plot[domain=-0.6283185307179586:0,variable=\t]({1*1.6180339887498951*cos(\t r)+0*1.6180339887498951*sin(\t r)},{0*1.6180339887498951*cos(\t r)+1*1.6180339887498951*sin(\t r)});
\draw [shift={(0,0)},line width=2pt]  plot[domain=0.6283185307179585:1.2566370614359172,variable=\t]({1*1.6180339887498947*cos(\t r)+0*1.6180339887498947*sin(\t r)},{0*1.6180339887498947*cos(\t r)+1*1.6180339887498947*sin(\t r)});
\draw [shift={(1,0)},line width=2pt]  plot[domain=1.884955592153876:2.5132741228718345,variable=\t]({1*1.6180339887498947*cos(\t r)+0*1.6180339887498947*sin(\t r)},{0*1.6180339887498947*cos(\t r)+1*1.6180339887498947*sin(\t r)});
\draw [shift={(1.3090169943749475,0.9510565162951532)},line width=2pt]  plot[domain=3.141592653589793:3.7699111843077513,variable=\t]({1*1.618033988749895*cos(\t r)+0*1.618033988749895*sin(\t r)},{0*1.618033988749895*cos(\t r)+1*1.618033988749895*sin(\t r)});
\draw [line width=.8pt,dashed] (0.5,1.5388417685876266)-- (0.5,-0.07919222016226816);
\draw [line width=.8pt,dashed] (0,0)-- (0.9510565162951535,1.3090169943749472);
\draw [line width=.8pt,dashed] (1,0)-- (0.04894348370484636,1.3090169943749475);
\draw [line width=.8pt,dotted] (0.5,0.6881909602355867) circle (0.8506508083520399cm);
\draw [shift={(0.5,0.6881909602355867)},line width=2pt]  plot[domain=3.455751918948772:4.71238898038469,variable=\t]({1*0.7673831803978549*cos(\t r)+0*0.7673831803978549*sin(\t r)},{0*0.7673831803978549*cos(\t r)+1*0.7673831803978549*sin(\t r)});
\draw [shift={(0.5,0.6881909602355867)},line width=2pt]  plot[domain=4.71238898038469:5.969026041820607,variable=\t]({1*0.7673831803978548*cos(\t r)+0*0.7673831803978548*sin(\t r)},{0*0.7673831803978548*cos(\t r)+1*0.7673831803978548*sin(\t r)});
\draw [shift={(0.5,0.6881909602355867)},line width=2pt]  plot[domain=-0.3141592653589793:0.9424777960769379,variable=\t]({1*0.7673831803978555*cos(\t r)+0*0.7673831803978555*sin(\t r)},{0*0.7673831803978555*cos(\t r)+1*0.7673831803978555*sin(\t r)});
\draw [shift={(0.5,0.6881909602355867)},line width=2pt]  plot[domain=0.9424777960769379:2.199114857512855,variable=\t]({1*0.7673831803978548*cos(\t r)+0*0.7673831803978548*sin(\t r)},{0*0.7673831803978548*cos(\t r)+1*0.7673831803978548*sin(\t r)});
\draw [shift={(0.5,0.6881909602355867)},line width=2pt]  plot[domain=2.199114857512855:3.455751918948772,variable=\t]({1*0.7673831803978551*cos(\t r)+0*0.7673831803978551*sin(\t r)},{0*0.7673831803978551*cos(\t r)+1*0.7673831803978551*sin(\t r)});
\begin{scriptsize}
\draw [fill=black] (0,0) circle (.2pt);
\draw[color=black] (-0.03,-0.03) node {$q_2$};
\draw [fill=black] (1,0) circle (.2pt);
\draw[color=black] (1.03,-0.03) node {$q_3$};
\draw [fill=black] (0.5,1.5388417685876266) circle (.2pt);
\draw[color=black] (0.5,1.58) node {$q_1$};
\draw [fill=black] (0.5,0.6881909602355867) circle (.2pt);
\draw[color=black] (0.518,0.75) node {$p$};
\draw [fill=black] (0.5,-0.07919222016226816) circle (.2pt);
\draw[color=black] (0.5,-0.12) node {$t_1=w_{23}=w_{32}$};
\draw [fill=black] (0.04894348370484636,1.3090169943749475) circle (.2pt);
\draw[color=black] (0.17,1.3) node {$t_3=w_{12}$};
\draw [fill=black] (0.9510565162951535,1.3090169943749472) circle (.2pt);
\draw[color=black] (0.83,1.3) node {$t_2=w_{13}$};
\draw [fill=black] (-0.22982477421267944,0.45105651629515375) circle (.2pt);
\draw[color=black] (-0.17,0.45) node {$w_{21}$};
\draw [fill=black] (1.2298247742126798,0.45105651629515314) circle (.2pt);
\draw[color=black] (1.17,0.45) node {$w_{31}$};
\end{scriptsize}
\end{tikzpicture}
\end{center}
Since for any $i=1,2,3$, $B(p,\varrho)\subset B(q_i,D)$ where $\{t_i\}=\partial B(p,\varrho)\cap \partial B(q_i,D)$, it follows from \eqref{pintt1t2t3}  that
\begin{equation}
\label{tiwijwik}
t_i\not\in {\rm int}\,B\left(z_{ij},D\right)
\mbox{ \ and \ }t_i\not\in {\rm int}\,B\left(z_{ik},D\right)\mbox{ for $\{i,j,k\}=\{1,2,3\}$.}
\end{equation}

We deduce from Lemma~\ref{CompleteProperties} that the shorter arc
$\arc{w_{ij}q_i}$ of $\partial B\left(z_{ij},D\right)$  is contained in $K$. Now let $\Gamma_i$ be the Jordan domain bounded by the shorter arc $\arc{w_{ij}q_i}$ of the circles $\partial B\left(z_{ij},D\right)$,
the shorter arc $\arc{w_{ik}q_i}$ of the circles $\partial B\left(z_{ik},D\right)$,
and the shorter arc $\arc{w_{ij}w_{ik}}$ of $B(p,\varrho)$, $\{i,j,k\}=\{1,2,3\}$.
Then $\Gamma_i\subseteq K$ for $i=1,2,3$,
and \eqref{tiwijwik} yields that
\begin{equation}
\label{GammaiGammajcap}
\mathrm{int}\Gamma_j\cap\mathrm{int}\Gamma_k=\emptyset
\end{equation}
for $i\neq j$. 
Now by the construction, $\Gamma_i$, $i=1,2,3$ are all congruent, and in turn congruent to $\widetilde{\Gamma}_k(\varrho)$, $j=1,2,3$
occurring at the definition of $\widetilde{Q}_D(\varrho)$, and hence \eqref{VQ} yields that
$$
V\left(K\right)\geq V\left(B\left(p,\varrho\right)\right)+
V\left(\Gamma_1\right)+V\left(\Gamma_2\right)+
V\left(\Gamma_3\right)=V\left(\widetilde{Q}_D\left(\varrho\right)\right).
$$
We conclude that
\begin{equation}
\label{VKVQVWVUD}
 V\left(K\right)\geq V\left(\widetilde{Q}_D\left(\varrho\right)\right)> V\left(\widetilde{W}_D\left(\varrho\right)\right)
 \geq V\left(U_D\right)   
\end{equation}
by \eqref{QDinWD} and Lemma~\ref{diskdomains}.
\hfill $\Box$\\

\section{Proof of Theorem~\ref{BL-stab}}

Let $\mathcal{M}^2$ be either $\R^2$, $S^2$ or $H^2$ and let $D>0$ where $D<\frac{\pi}2$ if $\mathcal{M}^2=S^2$.
In this section, we use the notation introduced in  Section~\ref{secReuleaux} and Section~\ref{sec-BL-hyp}.

The idea of the proof of Theorem~\ref{BL-stab}
based on the argument for Theorem~\ref{BL-hyp} runs as follows. If $K\subset \mathcal{M}^2$ is a convex body of constant width $D$, and
$
V_{\mathcal{M}^2}\left(K\right)\leq V_{\mathcal{M}^2}\left(U_D\right)+\varepsilon$ for small $\varepsilon>0$,
then using the notation of Theorem~\ref{BL-M2},
combining
$$
V_{\mathcal{M}^2}(\widetilde{Q}_D(r(K))=
V_{\mathcal{M}^2}\left(B\left(p,\varrho\right)\right)+
V_{\mathcal{M}^2}\left(\Gamma_1\right)+V\left(\Gamma_2\right)+
V_{\mathcal{M}^2}\left(\Gamma_3\right)\leq V_{\mathcal{M}^2}\left(K\right)\leq V_{\mathcal{M}^2}\left(U_D\right)+\varepsilon
$$
and Lemma~\ref{step5} below 
yields that $r(K)\leq r(U_D)+O(\varepsilon)$.
In order to verify Lemma~\ref{step5}, we provide a stability version of on Lemma~\ref{diskdomains}
about $\widetilde{W}_D(\varrho)\subset \widetilde{Q}_D(\varrho)$
if $\varrho$ is close to $r(U_D)$, and use that
$\widetilde{W}_D(\varrho)$ is strictly contained in
$\widetilde{Q}_D(\varrho)$ if $\varrho$ is away from $r(U_D)$.
Finally, we conclude via 
Proposition~\ref{step6} below that $K$ is close to a Reuleaux triangle of width $D$.

For the whole section, we fix $\eta_0=\frac{D}{4}-\frac{r\left(U_D\right)}{2}$.
We observe that if  $\varrho=r\left(U_D\right)+\eta$ for some $0<\eta\leq \eta_0$, then 
$$
r\left(U_D\right)<\varrho\leq \frac{D}{4}+\frac{r\left(U_D\right)}{2}<\frac{D}{2}.
$$

We prepare the proof of Theorem~\ref{BL-stab} with a series of lemmas.

\begin{lemma}
\label{step1}
For $\varrho=r\left(U_D\right)+\eta$ with $0<\eta\leq\eta_0$, there exists an explicit positive constant $\alpha_0$ depending on $D$ and $\mathcal{M}^2$ such that using the notation in the proof of Lemma~\ref{diskdomains}, we have
$$
 \angle (\widetilde{q}_1,v_1,u'_2)\geq
\alpha_0.  
$$
\end{lemma}
\proof We observe that
$$
\angle (\widetilde{q}_1,v_1,u'_2)=
\angle (\widetilde{p},v_1,\widetilde{p}')=
\angle (\widetilde{p},\widetilde{t}_2,\widetilde{p}')=
\angle (\widetilde{p},\widetilde{t}_2,v_1)-\angle (v_1,\widetilde{t}_2,\widetilde{p}')=
\angle (\widetilde{p},\widetilde{t}_2,v_1)-\angle (\widetilde{t}_2,v_1,\widetilde{p}).
$$
In the triangle $[\widetilde{p},\widetilde{t}_2,v_1]$,
we have $\angle (\widetilde{t}_2,\widetilde{p},v_1)=\frac{\pi}3$ and
\begin{eqnarray*}
d_{\mathcal{M}^2}(\widetilde{p},v_1)&=&
D-r\left(U_D\right)>\frac{D}{2}\\
d_{\mathcal{M}^2}(\widetilde{p},\widetilde{t}_2)&=&
\varrho\leq \frac{D}{4}+\frac{r\left(U_D\right)}{2}<
\frac{D}{2}.
\end{eqnarray*}
In particular, 
$\angle (\widetilde{q}_1,v_1,u'_2)=\angle (\widetilde{p},\widetilde{t}_2,v_1)-\angle (\widetilde{t}_2,v_1,\widetilde{p})$ is minimal if
$\widetilde{t}_2=\widetilde{t}_{2,0}$ where
$d_{\mathcal{M}^2}(\widetilde{p},\widetilde{t}_{2,0})=
\frac{D}{4}+\frac{r\left(U_D\right)}{2}$, and we set
$\alpha_0=\angle (\widetilde{p},\widetilde{t}_{2,0},v_1)-\angle (\widetilde{t}_{2,0},v_1,\widetilde{p})>0$.
\hfill $\Box$\\

Now let us give a lower bound for the geodesic distance
of $\widetilde{q}_1$ from $\gamma'_{12}$ inside
$\widetilde{W}_D(\varrho)$ ({\it cf.} Lemma~\ref{diskdomains}) if $\varrho$ is close to $r(U_D)$. We observe that 
$d_{\mathcal{M}^2}(v_2,v_1)=d_{\mathcal{M}^2}(v_2,u_2)=D$, and hence $d_{\mathcal{M}^2}(v'_2,\widetilde{t}_2)=d_{\mathcal{M}^2}(v_2,u'_2)=D$, and hence
$\gamma'_{12}\subset \partial B_{\mathcal{M}^2}\left(v_2',D\right)$
and $\sigma_{12}\subset\partial B_{\mathcal{M}^2}\left(\widetilde{z}_{12},D\right)$.

\begin{lemma}\label{step2}
If $0<\eta\leq\eta_0$ and $\varrho=r(U_D)+\eta$, then 
for the arc $\gamma'_{12}\subset \partial_{\mathcal{M}^2}B(v'_2,D)$
defined in the proof of Lemma~\ref{diskdomains},
the distance of $q_1$ from $B(v'_2,D)$, and hence
from $\gamma'_{12}$ is at least $\theta_1\eta$ for explicit $\theta_1>0$ depending on $D$ and $\mathcal{M}^2$.
\end{lemma}
\proof It is equivalent to prove that
$d_{\mathcal{M}^2}\left(\widetilde{q}_1,v_2'\right)
\geq D+\theta_1\eta$ for explicit $\theta_1>0$ depending on $D$ and $\mathcal{M}^2$; or in other words,
\begin{equation}
\label{q1vprime2dist}
\begin{array}{rcll}
 d_{\R^2}\left(\widetilde{q}_1,v_2'\right)^2  & \geq&
D^2+\widetilde{\theta}_1\eta &
\mbox{ if }\mathcal{M}^2=\R^2;\\
 \cos d_{S^2}\left(\widetilde{q}_1,v_2'\right)  & \leq&
\cos D-\widetilde{\theta}_1\eta &
\mbox{ if }\mathcal{M}^2=S^2;\\
\cosh d_{H^2}\left(\widetilde{q}_1,v_2'\right)  & \geq&
\cosh D+\widetilde{\theta}_1\eta &
\mbox{ if }\mathcal{M}^2=H^2.
\end{array}
\end{equation}
for explicit $\widetilde{\theta}_1>0$ depending on $D$ and $\mathcal{M}^2$.

Let $\alpha=\angle (\widetilde{q}_1,v_1,u'_2)$, and hence $\alpha\geq \alpha_0$ for the constant
$\alpha_0>0$ of Lemma~\ref{step1}. 
Using that
$$
d_{\mathcal{M}^2}\left(v_1,\widetilde{q}_1\right)=\varrho-r\left(U_D\right)=\eta
$$
and
$$
d_{\mathcal{M}^2}\left(v_1,v_2'\right)=d_{\mathcal{M}^2}\left(t_2,v_2\right)=R\left(U_D\right)+\varrho=D+\eta,
$$
we express the geodesic distance $d_{\mathcal{M}^2}\left(\widetilde{q}_1,v_2'\right)$ by applying the Law of Cosines for the triangle $\left[v_1,\widetilde{q}_1,v_2'\right]_{\mathcal{M}^2}$. The rest of the argument is divided into three cases depending on $\mathcal{M}^2$.

If $\mathcal{M}^2=\R^2$, then
\begin{eqnarray*}
d_{\R^2}\left(\widetilde{q}_1,v_2'\right)^2&=&\eta^2+\left(D+\eta\right)^2
-2\eta\left(D+\eta\right) \cos \alpha=D^2+
2\eta\left(D+\eta\right) (1-\cos \alpha)\\
&\geq &D^2+2D (1-\cos \alpha_0)\cdot \eta.
\end{eqnarray*}
If $\mathcal{M}^2=S^2$, then
\begin{eqnarray*}
\cos\left(d_{S^2}\left(\widetilde{q}_1,v_2'\right)\right)&=&
\cos\left(\eta\right)\cdot\cos\left(D+\eta\right)+\sin\left(\eta\right)\cdot\sin\left(D+\eta\right)\cdot\cos\alpha\\
&=&\cos\left((D+\eta)-\eta\right)-\sin\left(\eta\right)\cdot\sin\left(D+\eta\right)(1-\cos\alpha)\\
&\leq &
\cos D-\frac{(1-\cos\alpha_0)\sin D}4\cdot \eta.
\end{eqnarray*}
Finally, if $\mathcal{M}^2=H^2$, then
\begin{eqnarray*}
\cosh\left(d_{H^2}\left(\widetilde{q}_1,v_2'\right)\right)&=&
\cosh\left(\eta\right)\cdot\cosh\left(D+\eta\right)-
\sinh\left(\eta\right)\cdot\sinh\left(D+\eta\right)\cdot\cos\alpha\\
&=&\cosh\left((D+\eta)-\eta\right)+\sinh\left(\eta\right)\cdot\sinh\left(D+\eta\right)(1-\cos\alpha)\\
&\geq &
\cosh D+(1-\cos\alpha_0)\sinh D\cdot \eta,
\end{eqnarray*}
verifying \eqref{q1vprime2dist}, and in turn Lemma~\ref{step2}.
\hfill $\Box$\\

We will use the following area estimate in the proof of Proposition~\ref{step5} if
$\varrho=r\left(U_D\right)+\eta$ for $\eta\geq \eta_0$.

\begin{lemma}\label{step2half}
There exists an explicit positive constant $\theta_2$ depending on $D$ and $\mathcal{M}^2$ such that
if $r(U_D)+\eta_0\leq\varrho\leq D/2$, then
$$
V_{\mathcal{M}^2}\left(\widetilde{Q}_D\left(\varrho\right)\right)\geq V_{\mathcal{M}^2}\left(U_D\right)+\theta_2.
$$
\end{lemma}
\proof We observe that the largest angle
of the isosceles triangle whose sides are
$D-r(U_D)$, $D-r(U_D)$ and $D$ is $\frac{2\pi}{3}$. Therefore,
there exists $\vartheta_0>0$ such that
the largest angle
of the isosceles triangle whose sides are
$D-(r(U_D)+\eta_0)$, $D-(r(U_D)+\eta_0)$ 
and $D$ is $\frac{2\pi}{3}+\vartheta_0$; namely,
$$
\sin\left(\frac{\pi}{3}+\frac{\vartheta_0}2\right)=
\left\{\begin{array}{ll}
\frac{D/2}{D-(r(U_D)+\eta_0)} & 
\mbox{ \ if $\mathcal{M}^2=\R^2$}\\[1ex]
\frac{\sin(D/2)}{\sin\left(D-(r(U_D)+\eta_0)\right)} & 
\mbox{ \ if $\mathcal{M}^2=S^2$} \\[1ex]
\frac{\sinh(D/2)}{\sinh\left(D-(r(U_D)+\eta_0)\right)} & 
\mbox{ \ if $\mathcal{M}^2=H^2$}
\end{array} \right.
$$
Next we fix $a\in \mathcal{M}^2$, and choose
$b,c\in \partial_{\mathcal{M}^2}B_{\mathcal{M}^2}\left(a,\frac{D}2\right)$ such that
$\angle(b,a,c)=\vartheta_0$. In addition,
let $f\in \mathcal{M}^2$ satisfy that
$d_{\mathcal{M}^2}(b,f)=d_{\mathcal{M}^2}(c,f)=D$
and $a\in[b,c,f]_{\mathcal{M}^2}$. We claim that
$\theta_2=V_{\mathcal{M}^2}(\Xi_0)$ works in Lemma~\ref{step2half}
where $\Xi_0$ is the Jordan domain between the shorter arc of 
$\partial_{\mathcal{M}^2}B_{\mathcal{M}^2}\left(a,\frac{D}2\right)$ connecting $b$ and $c$
and the shorter arc of 
$\partial_{\mathcal{M}^2}B_{\mathcal{M}^2}\left(f,D\right)$ connecting $b$ and $c$.

To prove Lemma~\ref{step2half} with this $\theta_2$, we may assume that 
$r(U_D)+\eta_0<\varrho< D/2$. Let $\Theta$ be the part of $\widetilde{Q}_D(\varrho)$ cut off by the circular arc $\sigma_{12}$ of radius $D$ connecting $\widetilde{q}_1$ and $\widetilde{t}_2$ on the boundary of $\widetilde{W}_D(\varrho)$; namely, $\Theta$ is a Jordan domain bounded by $\sigma_{12}$, by the arc $\arc{\widetilde{t}_2\widetilde{w}_{12}}\subset \partial_{\mathcal{M}^2}B_{\mathcal{M}^2}(\widetilde{p},\varrho)$ and the arc
$\arc{\widetilde{w}_{12}\widetilde{q}_1}\subset \partial_{\mathcal{M}^2}B_{\mathcal{M}^2}(\widetilde{z}_{ij},D)$, and satisfies
$$
V_{\mathcal{M}^2}\left(\widetilde{Q}_D(\varrho)\right)
=V_{\mathcal{M}^2}\left(\widetilde{Q}_D(\varrho)\right)
+6V_{\mathcal{M}^2}(\Theta).
$$
Since $B_{\mathcal{M}^2}\left(\widetilde{p},\varrho\right)$ and 
$B_{\mathcal{M}^2}\left(\widetilde{z}_{12},D\right)$
are touching each other in $\widetilde{w}_{12}$,
it follows that $\widetilde{p}\in[\widetilde{z}_{12},\widetilde{w}_{12}]_{\mathcal{M}^2}$, and hence
$$
d_{\mathcal{M}^2}(\widetilde{z}_{12},\widetilde{p})
=D-\varrho=d_{\mathcal{M}^2}(\widetilde{q}_{1},\widetilde{p}).
$$
In particular, the sides of the
triangle $[\widetilde{z}_{12},\widetilde{p},\widetilde{q}_{1}]_{\mathcal{M}^2}$ are of length $D-\varrho$, $D-\varrho$ and $D$, and $\varrho>r(U_D)+\eta_0$
yields that the angle
$\angle(\widetilde{z}_{12},\widetilde{p},\widetilde{q}_{1})$ is of the form $\frac{2\pi}3+\vartheta$ for
$\vartheta>\vartheta_0$. As
$\angle(\widetilde{t}_{2},\widetilde{p},\widetilde{q}_{1})=\frac{\pi}3$, we deduce that
$\angle(\widetilde{w}_{12},\widetilde{p},\widetilde{t}_{2})=\vartheta$.

Let $\xi(D/2)$ be the shorter circular arc 
of radius $D/2$ connecting $\widetilde{w}_{12}$
and $\widetilde{t}_{2}$ and being separated by
the segment $[\widetilde{w}_{12},\widetilde{t}_{2}]_{\mathcal{M}^2}$ from $\widetilde{p}$, let $\xi(D)$ be the shorter circular arc 
of radius $D$ connecting $\widetilde{w}_{12}$
and $\widetilde{t}_{2}$ and being separated by
the segment $[\widetilde{w}_{12},\widetilde{t}_{2}]_{\mathcal{M}^2}$ from $\widetilde{p}$, and let
$\Xi$ be the Jordan domain bounded by $\xi(D/2)$ and $\xi(D)$. Since $\varrho< D/2$, it follows that
$$
\Xi\subset \Theta\subset \widetilde{Q}_D(\varrho)\backslash \widetilde{W}_D(\varrho).
$$
On the other hand, $\angle(\widetilde{w}_{12},\widetilde{p},\widetilde{t}_{2})=\vartheta>\vartheta_0$ implies that
$\Xi$ contains a congruent copy of $\Xi_0$; therefore,
$$
V_{\mathcal{M}^2}\left(\widetilde{Q}_D(\varrho)\right)
\geq V_{\mathcal{M}^2}\left(\widetilde{Q}_D(\varrho)\right)
+6V_{\mathcal{M}^2}(\Xi_0)>V_{\mathcal{M}^2}\left(\widetilde{Q}_D(\varrho)\right)+\theta_2.
\mbox{ \ }\Box
$$

One of the main auxiliary statement towards proving
Theorem~\ref{BL-stab} is the following statement.

\begin{prop}\label{step5}
There exists an explicit positive constant $\theta_3$ depending on $D$ and $\mathcal{M}^2$ such that
if $r(U_D)\leq \varrho\leq D/2$ and 
$\varrho=r\left(U_D\right)+\eta$, then
$$
V_{\mathcal{M}^2}\left(\widetilde{Q}_D\left(\varrho\right)\right)\geq V_{\mathcal{M}^2}\left(U_D\right)+\theta_3\eta.
$$
\end{prop}
\proof If $\eta\geq \eta_0$, then 
Lemma~\ref{step2half} yields that any $\theta_3>0$
with
\begin{equation}
\label{theta2theta3eta0}
\theta_3\leq \frac{\theta_2}{\frac{D}2-r(U_D)}
\end{equation}
works in Proposition~\ref{step5}. Therefore, we assume that 
$$
\varrho=r\left(U_D\right)+\eta
\mbox{ \ for }\eta\leq \eta_0.
$$
Since $\widetilde{W}_D\left(\varrho\right)\subset 
\widetilde{Q}_D\left(\varrho\right)$,
and we have seen in the proof of Lemma~\ref{diskdomains} that
$$
\widetilde{W}_D\left(\varrho\right)\setminus U_D\supsetneq\Omega_+\setminus\widetilde{\Omega},
$$
it is sufficient to prove that if 
$\eta\leq \eta_0$, then
\begin{equation}\label{Omega_area}
V_{\mathcal{M}^2}\left(\Omega_+\setminus\widetilde{\Omega}\right)\geq\widetilde{\theta}_3\eta.
\end{equation}
for explicit $\widetilde{\theta}_3>0$
depending on $D$ and $\mathcal{M}^2$.
We recall that the closure of $\Omega_+\setminus\widetilde{\Omega}$ is the Jordan domain bounded by
the arc $\arc{\widetilde{t}_{12},e}$ of $\gamma'_{12}$ where $e=\gamma_{12}\cap \gamma'_{12}$,
the arc $\arc{\widetilde{t}_{12},d}$ of $\sigma_{12}$ where $d=\gamma_{12}\cap \sigma_{12}$
and the arc $\arc{d,e}$ of $\gamma_{12}$.

To verify \eqref{Omega_area}, first
we claim that the 
angle 
$\angle(v'_2,\widetilde{t}_2,\widetilde{z}_{12})$ of the tangents of the circles $\partial B_{\mathcal{M}^2}\left(v_2',D\right)$ and $\partial B_{\mathcal{M}^2}\left(\widetilde{z}_{12},D\right)$ at $\widetilde{t}_2$ satisfies
\begin{equation}
\label{anglecircles}
\angle(v'_2,\widetilde{t}_2,\widetilde{z}_{12})\geq \theta_4\eta
\end{equation}
for some explicit positive constant $\theta_4$ depending only on $D$ and $\mathcal{M}^2$. It follows from
Lemma~\ref{step2} that 
$d_{\mathcal{M}^2}(v'_2,q_1)\geq D+\theta_2\eta$ where
$q_1\in \sigma_{12}\subset \partial B_{\mathcal{M}^2}\left(\widetilde{z}_{12},D\right)$; therefore, the triangle inequality yields that 
$d_{\mathcal{M}^2}(v'_2,\widetilde{z}_{12})\geq \theta_2\eta$. 
Writing $m$ to denote the midpoint of the segment $[v'_2,\widetilde{z}_{12}]_{\mathcal{M}^2}$, 
we have $\angle(v'_2,m,\widetilde{t}_2,m)=\frac{\pi}2$, and hence
the sine of angle $\angle(v'_2,\widetilde{t}_2,m)=\frac12\angle(v'_2,\widetilde{t}_2,\widetilde{z}_{12})$
can be expressed in terms of $d_{\mathcal{M}^2}(v'_2,m)\geq \frac{\theta_2}2\cdot \eta$
and $d_{\mathcal{M}^2}(v'_2,\widetilde{t}_2)=D$ using the Law of Sines for the triangle
$[v'_2,\widetilde{t}_2,m]_{\mathcal{M}^2}$. In turn, we conclude \eqref{anglecircles}.

As $\angle(\widetilde{t}_2,\widetilde{p},v_1)=\frac{\pi}3$,
$d_{\mathcal{M}^2}(\widetilde{p},v_1)=D-r(U_D)>\frac{D}2$,
we deduce that
$d_{\mathcal{M}^2}(\widetilde{t}_2,v_1)>\frac{D}4$. Thus the property
$d_{\mathcal{M}^2}(\widetilde{t}_2,e)=d_{\mathcal{M}^2}(v_1,e)$  of the point $e=\gamma_{12}\cap \gamma'_{12}$ yields that
$d_{\mathcal{M}^2}(\widetilde{t}_2,v_1)>\frac{D}8$. Therefore if $a_\gamma\in \gamma'_{12}$ has  the property that 
$d_{\mathcal{M}^2}(\widetilde{t}_2,a_\gamma)=\frac{D}{16}$ and $a_\sigma\in \sigma_{12}$ has  the property
that $d_{\mathcal{M}^2}(\widetilde{t}_2,a_\sigma)=\frac{D}{16}$, then
$$
\Psi\subset {\rm cl}\,\left(\Omega_+\setminus\widetilde{\Omega}\right)
$$
for the Jordan domain $\Psi$ bounded by the arc $\arc{\widetilde{t}_{12},a_\gamma}$ of $\gamma'_{12}$,
the arc $\arc{\widetilde{t}_{12},a_\sigma}$ of $\sigma'_{12}$ and the shorter 
$\arc{a_\gamma,a_\sigma}$ of the circle 
$\partial B_{\mathcal{M}^2}\left(\widetilde{t}_{12},\frac{D}{16}\right)$. Here  
 $\angle(a_\gamma,\widetilde{t}_{12},a_\sigma)$ equals the angle of the arcs
$\gamma'_{12}$ and $\sigma_{12}$, which is at least $\theta_4\eta$ according to \eqref{anglecircles}.
We observe that $V_{\mathcal{M}^2}(\Psi)=V_{\mathcal{M}^2}\left(\widetilde{\Psi}\right)$ for the circular sector
$\widetilde{\Psi}$ of $B_{\mathcal{M}^2}\left(\widetilde{t}_{12},\frac{D}{16}\right)$  bounded by the segments
$[\widetilde{t}_{12},a_\sigma]_{\mathcal{M}^2}$
and $[\widetilde{t}_{12},a_\gamma]_{\mathcal{M}^2}$
where
$$
V_{\mathcal{M}^2}\left(\widetilde{\Psi}\right)=
\frac{\angle(a_\gamma,\widetilde{t}_{12},a_\sigma)}{2\pi}\cdot V_{\mathcal{M}^2}\left(B_{\mathcal{M}^2}\left(\widetilde{t}_{12},\frac{D}{16}\right)\right)\geq
\frac{\theta_4V_{\mathcal{M}^2}\left(B_{\mathcal{M}^2}\left(\widetilde{t}_{12},\frac{D}{16}\right)\right)}{2\pi}
\cdot \eta.
$$
We deduce that
$$
V_{\mathcal{M}^2}\left(\Omega_+\setminus\widetilde{\Omega}\right)\geq
V_{\mathcal{M}^2}(\Psi)=V_{\mathcal{M}^2}\left(\widetilde{\Psi}\right)\geq
\frac{\theta_4V_{\mathcal{M}^2}\left(B_{\mathcal{M}^2}\left(\widetilde{t}_{12},\frac{D}{16}\right)\right)}{2\pi}
\cdot\eta.
$$
We conclude the claim \eqref{Omega_area}, which in turn yields Proposition~\ref{step5} by \eqref{theta2theta3eta0}.
\hfill $\Box$\\

\begin{lemma}\label{step6}
There exist explicit constants $\theta_4>0$
and $\widetilde{\eta}_0\in(0,\eta_0]$ depending on $D$ and $\mathcal{M}^2$ such that
if $\varrho=r(U_D)+\eta$ for $0<\eta\leq \widetilde{\eta}_0$ and $\{i,j,k\}=\{1,2,3\}$, then  
\begin{equation}\label{step6eq}
\angle(\widetilde{w}_{ij},\widetilde{p},\widetilde{w}_{ik})\geq \frac{2\pi}3-\theta_4\eta.
\end{equation}
\end{lemma}

\proof We may assume that $i=1$.
The proof of Lemma~\ref{step6} is rather technical; therefore, we sketch the main idea. 
Let $\widetilde{v}_1$ be the intersection of $\partial B(\widetilde{q}_2,D)$ and the half-line $\widetilde{p}v_1$, and hence $v_1\in [\widetilde{p},\widetilde{v}_1]_{\mathcal{M}^2}$.
Our first claim is that there exists an explicit positive constant $\theta_5$ depending on $D$ and $\mathcal{M}^2$
such that if $\eta\leq \eta_0$, then
\begin{equation}
\label{v1tildev1}
d_{\mathcal{M}^2}(\widetilde{v}_1,v_1)\leq \theta_5\eta,
\end{equation}
which in turn yields $d_{\mathcal{M}^2}(\widetilde{v}_1,\widetilde{q}_1)\leq (\theta_5+1)\eta$
as $d_{\mathcal{M}^2}(v_1,\widetilde{q}_1)=\eta$.

Next the circle $\partial B_{\mathcal{M}^2}(\widetilde{q}_2,D)$ passing through $\widetilde{v}_1$ touches the circle
$B_{\mathcal{M}^2}(\widetilde{p},\varrho)$ at $\widetilde{t}_2$. Let 
$\widetilde{a}$ be the point of the shorter arc $\arc{\widetilde{t}_2,\widetilde{v}_1}$ 
of  $\partial B_{\mathcal{M}^2}(\widetilde{q}_2,D)$ satisfying
$d_{\mathcal{M}^2}(\widetilde{p},\widetilde{a})=d_{\mathcal{M}^2}(\widetilde{p},\widetilde{q}_1)$.
Since the circle $\partial B_{\mathcal{M}^2}(\widetilde{z}_{12},D)$ passing through $\widetilde{q}_1$ touches the circle
$B_{\mathcal{M}^2}(\widetilde{p},\varrho)$ at $\widetilde{w}_{12}$, we deduce that the triangles
$[\widetilde{t}_2,\widetilde{p},\widetilde{a}]_{\mathcal{M}^2}$ and $[\widetilde{w}_{12},\widetilde{p},\widetilde{q}_1]_{\mathcal{M}^2}$ are congruent where $\angle(\widetilde{t}_2,\widetilde{p},\widetilde{a})=\angle(\widetilde{w}_{12},\widetilde{p},\widetilde{q}_1)$.
Since $\angle(\widetilde{w}_{12},\widetilde{p},\widetilde{w}_{13})=2\cdot \angle(\widetilde{w}_{12},\widetilde{p},\widetilde{q}_1)$,
Lemma~\ref{step6} is equivalent proving that
there exist explicit constants $\theta_4>0$
and $\widetilde{\eta}_0\in(0,\eta_0]$ depending on $D$ and $\mathcal{M}^2$ such that
if $\varrho=r(U_D)+\eta$ for $\eta\in[0,\widetilde{\eta}_0]$, then
\begin{equation}
\label{step6eq0}
\angle(\widetilde{t}_2,\widetilde{p},\widetilde{a})\geq \frac{\pi}3-\frac{\theta_4}2\cdot\eta,
\end{equation}
which statement we verify using \eqref{v1tildev1}; more precisely, using $d_{\mathcal{M}^2}(\widetilde{v}_1,\widetilde{q}_1)\leq (\theta_5+1)\eta$.

To prove \eqref{v1tildev1}, we observe that $\angle(\widetilde{v}_1,\widetilde{p},\widetilde{q}_2)=\frac{2\pi}3$, $d_{\mathcal{M}^2}(\widetilde{v}_1,\widetilde{q}_2)=D$
 and $d_{\mathcal{M}^2}(\widetilde{q}_2,\widetilde{p})=R(U_D)-\eta$
in the triangle $[\widetilde{v}_1,\widetilde{p},\widetilde{q}_2]_{\mathcal{M}^2}$, and
$\angle(v_1,\widetilde{p},v_2)=\frac{2\pi}3$, 
$d_{\mathcal{M}^2}(v_1,\widetilde{p})=R(U_D)$,
$d_{\mathcal{M}^2}(v_1,v_2)=D$ and $d_{\mathcal{M}^2}(v_2,\widetilde{p})=R(U_D)$
in the triangle $[v_1,\widetilde{p},v_2]_{\mathcal{M}^2}$. For $s\in\left[\frac{D}2,R(U_D)\right]$, let $f(s)\in[R(U_D),D]$ be defined in a way such that there exists a triangle with side lengths $s,f(s),D$ and having angle $\frac{2\pi}3$ opposite to the side of length $D$. In particular,
$f(R(U_D)-\eta)=d_{\mathcal{M}^2}(\widetilde{v}_1,\widetilde{p})$ and 
$f(R(U_D))=d_{\mathcal{M}^2}(v_1,\widetilde{p})$, and the Law of Cosines yields
$$
\begin{array}{rcll}
D^2&=&s^2+f(s)^2+s\cdot f(s)&\mbox{ if $\mathcal{M}^2=\R^2$}\\
\cos D&=&\cos s\cdot \cos f(s)-\frac12\cdot \sin s\cdot \sin f(s)&\mbox{ if $\mathcal{M}^2=S^2$}\\
\cosh D&=&\cosh s\cdot \cosh f(s)+\frac12\cdot\sinh s\cdot \sinh f(s)&\mbox{ if $\mathcal{M}^2=H^2$}.
\end{array}
$$
It follows that
$$
f'(s)=\left\{
\begin{array}{ll}
-\frac{2s+f(s)}{s+2f(s)}&\mbox{ if $\mathcal{M}^2=\R^2$}\\[1ex]
-\frac{\sin s\cdot \cos f(s)+\frac12\cdot \cos s\cdot \sin f(s)}{\cos s\cdot \sin f(s)+\frac12\cdot \sin s\cdot \cos f(s)}&\mbox{ if $\mathcal{M}^2=S^2$}\\[1ex]
-\frac{\sinh s\cdot \cosh f(s)+\frac12\cdot\cosh s\cdot \sinh f(s)}{\cosh s\cdot \sinh f(s)+\frac12\cdot\sinh s\cdot \cosh f(s)}&\mbox{ if $\mathcal{M}^2=H^2$};
\end{array}\right.
$$
therefore, $-\theta_5\leq f'(s)\leq 0$ for an explicit positive constant $\theta_5$ depending on $D$ and $\mathcal{M}^2$.
As $d_{\mathcal{M}^2}(\widetilde{v}_1,v_1)=f(R(U_D)-\eta)-f(R(U_D))$, we conclude \eqref{v1tildev1}.

The estimate that we really need is not \eqref{v1tildev1}, but its corollary 
\begin{equation}
\label{tildeq1tildev1}
d_{\mathcal{M}^2}(\widetilde{v}_1,\widetilde{q}_1)\leq (\theta_5+1)\eta,
\end{equation}
which follows from $d_{\mathcal{M}^2}(v_1,\widetilde{q}_1)=\eta$.

Next, we prove \eqref{step6eq0}.
For $\beta\in\left[\frac{2\pi}3,\pi\right]$, let
$a(\beta)$ be the point of the shorter arc of  $\partial B_{\mathcal{M}^2}(\widetilde{q}_2,D)$ satisfying
$\angle(a(\beta),\widetilde{p},\widetilde{q}_2)=\beta$, and hence $a\left(\frac{2\pi}3\right)=\widetilde{v}_1$
and $a(\pi)=\widetilde{t}_2$. Let $\widetilde{\beta}\in\left(\frac{2\pi}3,\pi\right)$ satisfy that $\widetilde{a}=a(\widetilde{\beta})$.

We observe that $\angle(a(\beta),\widetilde{p},\widetilde{q}_2)=\beta$, 
$d_{\mathcal{M}^2}(a(\beta),\widetilde{q}_2)=D$
 and $d_{\mathcal{M}^2}(\widetilde{q}_2,\widetilde{p})=R(U_D)-\eta$
in the triangle $[a(\beta),\widetilde{p},\widetilde{q}_2]_{\mathcal{M}^2}$.
We set 
$g(\beta)=d_{\mathcal{M}^2}(a(\beta),\widetilde{p})$
and $R_\eta=R(U_D)-\eta$,
and hence 
 the Law of Cosines in the triangle $[a(\beta),\widetilde{p},\widetilde{q}_2]_{\mathcal{M}^2}$ yields
$$
\begin{array}{rcll}
D^2&=&R_\eta^2+g(\beta)^2-R_\eta\cdot g(\beta)\cdot\cos\beta&\mbox{ if $\mathcal{M}^2=\R^2$}\\
\cos D&=&\cos R_\eta\cdot \cos g(\beta)+ \sin R_\eta\cdot \sin g(\beta)\cdot\cos\beta&\mbox{ if $\mathcal{M}^2=S^2$}\\
\cosh D&=&\cosh R_\eta\cdot \cosh g(\beta)-\sinh R_\eta\cdot \sinh g(\beta)\cdot\cos\beta&\mbox{ if $\mathcal{M}^2=H^2$}.
\end{array}
$$
It follows that 
$$
g'(\beta)=\left\{\begin{array}{ll}
-\frac{R_\eta g(\beta)\sin\beta}{2g(\beta)-R_\beta\cos\beta}&\mbox{ if $\mathcal{M}^2=\R^2$}\\[1ex]
-\frac{\sin R_\eta\cdot \sin g(\beta)\cdot \sin\beta}
{\cos R_\eta\cdot \sin g(\beta)-\sin R_\eta\cdot \cos g(\beta)\cdot \cos\beta}&\mbox{ if $\mathcal{M}^2=S^2$}\\[1ex]
-\frac{\sinh R_\eta\cdot\sinh g(\beta)\cdot\sin\beta}{\cosh R_\eta\cdot\sinh g(\beta)-\sinh R_\eta\cdot \cosh g(\beta)\cdot \cos\beta}&\mbox{ if $\mathcal{M}^2=H^2$};
\end{array}\right.
$$
where $\cos\beta<0$ (as $\beta\geq \frac{2\pi}3$), $r(U_D)\leq g(U_D)\leq D$ and $\frac{D}2\leq R_\eta\leq D$.
Therefore, $g'(\beta)<0$ if $\beta\geq \frac{2\pi}3$, and
there exits an explicit positive constant $\theta_6$ depending on $D$ and $\mathcal{M}^2$
such that $g'(\beta)\leq -\theta_6$ if $\beta\in\left[\frac{2\pi}3,\frac{5\pi}6\right]$.
We observe that 
$g\left(\frac{2\pi}3\right)=d_{\mathcal{M}^2}(\widetilde{v}_1,\widetilde{p})$, 
$$
g(\widetilde{\beta})=d_{\mathcal{M}^2}(a(\widetilde{\beta}),\widetilde{p})=d_{\mathcal{M}^2}(\widetilde{q}_1,\widetilde{p})
\geq g\left(\frac{2\pi}3\right)-(\theta_5+1)\eta
$$
by \eqref{tildeq1tildev1}, and
$(\theta_5+1)\eta\leq \theta_6\cdot\frac{\pi}6$ if $\eta\leq \widetilde{\eta}_0$
for
$$
\widetilde{\eta}_0=\min\left\{\eta_0,\frac{\theta_6\cdot\frac{\pi}6}{\theta_5+1}\right\}.
$$
Therefore, if $0\leq \eta\leq \widetilde{\eta}_0$, then
$$
\angle(\widetilde{t}_2,\widetilde{p},\widetilde{a})=\pi-\widetilde{\beta}=
\frac{\pi}3-\left(\widetilde{\beta}-\frac{2\pi}3\right)\geq \frac{\pi}3-
\frac{\theta_5+1}{\theta_6}\cdot\eta,
$$
and hence we may choose $\theta_4=2\cdot \frac{\theta_5+1}{\theta_6}$ in
\eqref{step6eq0}.
In turn, we conclude Lemma~\ref{step6}.
\hfill $\Box$

Now we are all set to prove our main result, Theorem~\ref{BL-stab}. For the reader's convenience, we also restate Theorem~\ref{BL-stab} in an equivalent form that suits better the actual argument.

\bigskip\noindent\textbf{Theorem~\ref{BL-stab}. }\emph{Let $\mathcal{M}^2$ be either $\R^2$, $S^2$ or $H^2$ and let $D>0$ where $D<\frac{\pi}2$ if $\mathcal{M}^2=S^2$.
If $K\subset \mathcal{M}^2$ is a body of constant width $D$, $\varepsilon\geq 0$ and  
$$
V_{\mathcal{M}^2}\left(K\right)\leq V_{\mathcal{M}^2}\left(U_D\right)+\varepsilon,
$$
then there exists a Reuleaux triangle $U\subset \mathcal{M}^2$ of width $D$ such that
$\delta_H(K,U)\leq \theta\varepsilon$ where $\theta>0$ is an explicitly calculable constant depending on $D$ and $\mathcal{M}^2$.}
\proof
Let $K$ be convex body in $\mathcal{M}^2$ of constant width $D$ with $V(K)<V(U_D)+\varepsilon$. 
During the argument, we use the notation of Theorem~\ref{BL-M2}.
First we assume that $0<\varepsilon\leq \varepsilon_0$ for $\varepsilon_0=\theta_3\widetilde{\eta}_0$
where $\theta_3$ comes from Proposition~\ref{step5} and $\widetilde{\eta}_0$ comes from Lemma~\ref{step6}.

Let $\varrho=r(K)$ and let $Q$ be the union of $B(p,\varrho),\Gamma_1,\Gamma_2,\Gamma_3$, and hence $V(Q)=V(\widetilde{Q}_D(\varrho))$. It follows from Proposition~\ref{step5} that $\varrho=r(U_D)+\eta$ where $\eta\leq\theta_3^{-1}\varepsilon\leq \widetilde{\eta}_0$. We deduce from Lemma~\ref{step6} that 
if $\{i,j,k\}=\{1,2,3\}$, then
$$
 \angle(w_{ij},p,w_{ik})\geq \frac{2\pi}3-\theta_7\varepsilon
$$
for $\theta_7=\theta_4\theta_3^{-1}$ where $\theta_4$ comes from Lemma~\ref{step6}.
As the interiors of the triangles $[w_{12},p,w_{13}]_{\mathcal{M}^2}$, $[w_{22},p,w_{23}]_{\mathcal{M}^2}$
and $[w_{31},p,w_{32}]_{\mathcal{M}^2}$ are pairwise disjoint by \eqref{GammaiGammajcap}.
Since the segment $[p,q_i]_{\mathcal{M}^2}$ halves the angle $\angle(w_{ij},p,w_{ik})$
when $\{i,j,k\}=\{1,2,3\}$, we deduce that
if $i\neq j$, then
\begin{equation}\label{angles_in_Q}
\frac{2\pi}3-\theta_7\varepsilon\leq \angle(q_i,p,q_j)\leq \frac{2\pi}3+2\theta_7\varepsilon.
\end{equation}

We may assume that that triangles $[q_1,q_2,q_3]_{\mathcal{M}^2}$ and 
$[\widetilde{q}_1,\widetilde{q}_2,\widetilde{q}_3]_{\mathcal{M}^2}$ have the same orientation, $p=\widetilde{p}$
and $q_1=\widetilde{q}_1$. First $\varrho\leq r(U_D)+\theta_3^{-1}\varepsilon$ yields that
\begin{equation}\label{UDtildeQD}
\delta_H\left(\widetilde{Q}_D\left(\varrho\right),U_D\right)<\theta_{9}\varepsilon
\end{equation}
for an explicit $\theta_9>0$ depending on $D$ and $\mathcal{M}^2$.
We deduce from \eqref{angles_in_Q} that
$\angle(q_2,p,\widetilde{q}_2)\leq 2\theta_7\varepsilon$ and 
$\angle(q_3,p,\widetilde{q}_3)\leq 2\theta_7\varepsilon$; therefore,
\begin{equation}\label{QDtildeQD}
\delta_H\left(\widetilde{Q}_D\left(\varrho\right),U_D\right)<\theta_{10}\varepsilon
\end{equation}
for an explicit $\theta_{10}>0$ depending on $D$ and $\mathcal{M}^2$. Let
$$
M=\bigcap_{i=1}^3B_{\mathcal{M}^2}(\widetilde{q}_i,D+\theta_{10}\varepsilon), 
$$
and hence $\widetilde{Q}_D\subset M$ and $K\subset M$ where we use \eqref{QDtildeQD} and that
the diameter of $K$ is at most $D$ to verify $K\subset M$. Since 
\begin{equation}\label{MtildeQD}
\delta_H\left(\widetilde{Q}_D\left(\varrho\right),M\right)<\theta_{11}\varepsilon
\end{equation}
for an explicit $\theta_{11}>0$ depending on $D$ and $\mathcal{M}^2$, we deduce from
\eqref{UDtildeQD}, \eqref{QDtildeQD}, \eqref{MtildeQD} and
$Q_D\subset K\subset M$ that
\begin{equation}\label{KUDdist}
\delta_H\left(K,U_D\right)<\theta_{12}\varepsilon
\end{equation}
for an explicit $\theta_{12}>0$ depending on $D$ and $\mathcal{M}^2$.

Finally, if $\varepsilon\geq \varepsilon_0$, then we can position $K$ and $U_D$ is a way such that
$\delta_H\left(K,U_D\right)\leq 2D$; therefore,
we can choose $\theta=\max\left\{\theta_{12},\frac{2D}{\varepsilon_0}\right\}$ in Theorem~\ref{BL-stab}.
\hfill $\Box$

\noindent{\bf Acknowledgement: } 
K\'aroly J. B\"or\"oczky and \'Ad\'am Sagmeister are supported by NKFIH project  
 K 132002.

\end{document}